\RequirePackage{setspace}
\documentclass[aap,preprint]{imsart}
\usepackage{amsmath,amssymb,color}
\usepackage[T1]{fontenc}

\usepackage[francais]{babel}
\usepackage{paralist}
\usepackage{mathrsfs}
\usepackage{t1enc}
\usepackage[utf8]{inputenc}  
\usepackage[T1]{fontenc}    
\usepackage{hyperref}
\usepackage{cases}
\usepackage{setspace}
\DeclareMathAlphabet{\mathpzc}{OT1}{pzc}{m}{it}
\newcommand{\bS}{\mathbb{S}} 
\DeclareMathOperator{\sech}{sech}
\DeclareMathOperator{\ch}{cosh}
\DeclareMathOperator{\sh}{sinh}
\newcommand{\Z}{\mathbb{Z}}
\newcommand{\N}{\mathbb{N}}
\newcommand{\R}{\mathbb{R}}
\newcommand{\T}{\mathbb{T}}
\newcommand{\dint}{\displaystyle\int}
\newcommand{\eps}{\varepsilon}
\newcommand{\norm}[2]{\left\|#1\right\|_{#2}}
\newcommand{\abs}[1]{\left|#1\right|}
\newcommand{\BigOh}[1]{\mathcal{O}(#1)}
\def \H{I\!\!H}
\usepackage{color}
\newcommand{\cf}{\mathcal{F}}
\DeclareMathAlphabet{\mathpzc}{OT1}{pzc}{m}{it}

\startlocaldefs

\newtheorem{theorem}{Theorem}[part]
\newtheorem{definition}{Definition}[part]
\newtheorem{Proposition}{Proposition}[part]

\newtheorem{lemma}{Lemma}[part]
\newtheorem{corollary}{Corollary}[part]
\newtheorem{remark}{Remark}[part]

\newcommand{\vertiii}[1]{{\left\vert\kern-0.25ex\left\vert\kern-0.25ex\left\vert #1 
		\right\vert\kern-0.25ex\right\vert\kern-0.25ex\right\vert}}

\newcommand{\vvertiii}[1]{{\vert\kern-0.25ex\vert\kern-0.25ex\vert #1 
		\vert\kern-0.25ex\vert\kern-0.25ex\vert}}

\newcommand{\Q}{\mathbb{Q}}
\newcommand{\ce}{\mathcal{E}}
\newcommand{\ct}{\mathcal{T}}
\newcommand{\stopt}{\mathcal{T}_{t,T}}
\newcommand{\stops}{\mathcal{T}_{S}^p}
\newcommand{\stopo}{\mathcal{T}_0^p}
\newcommand{\stopsp}{\mathcal{T}_{S^+}^p}
\newcommand{\stopop}{\mathcal{T}_{0^+}^p}
\newcommand{\p}{\mathbb{P}}
\DeclareMathOperator*{\esssup}{ess\,sup} \DeclareMathOperator{\e}{e}
\DeclareMathOperator*{\essinf}{ess\,inf} \DeclareMathOperator{\f}{f}

\newcommand{\dbarY}{\overline{\overline{Y}}}
\newcommand{\barY}{\overline{Y}}
\newcommand{\limsupn}{\limsup_{n\rightarrow\infty}}
\newcommand{\limn}{\lim_{n\rightarrow\infty}}
\newcommand{\limup}{\lim_n\!\smash{\uparrow}}

\newcommand{\dbarA}{\overline{\overline{A}}}
\def\argmax{\mathop{\mathrm{argmax}}}

\def\RB{\mathbb{R}}
\def\DB{\mathbb{D}}
\def\EB{\mathbb{E}}

\def\FB{\mathbb{F}}

\def\Fc{\mathcal{F}}
\def\PC{\mathcal{P}}
\def\EC{\mathcal{E}}
\def\AC{\mathcal{A}}
\def\UC{\mathcal{U}}
\def\CC{\mathcal C}
\def\RC{\mathcal R}
\def\HC{\mathcal H}
\def\VC{\mathcal V}
\def\R{{\bf R}}

\def\coz{\mathrm{coz}}

\def\tY{\tilde{Y}}
\def\tN{\tilde{N}}
\def\tl{\tilde{\lambda}}
\def\tH{\tilde{H}}
\def\tp{\tilde{p}}
\def\tr{\tilde{r}}
\def\tq{\tilde{q}}
\def\tZ{\tilde{Z}}
\def\tK{\tilde{K}}
\def\hX{\hat{X}}
\def\hZ{\hat{Z}}
\def\hA{\hat{A}}
\def\hY{\hat{Y}}
\def\hu{\hat{u}}
\def\hb{\hat{b}}
\def\hg{\hat{g}}
\def\hc{\hat{c}}
\def\hQ{\hat{Q}}
\def\hH{\hat{H}}
\def\hp{\hat{p}}
\def\hth{\hat{\t}}
\def\hq{\hat{q}}
\def\hr{\hat{r}}
\def\hxi{\hat{\xi}}

\def\vZ{\stackrel{\vee}{Z}}  
\def\vY{\stackrel{\vee}{Y}}
\def\vK{\stackrel{\vee}{K}}

\def\vp{\stackrel{\vee}{\pi}}
\def\vl{\stackrel{\vee}{l}}

\def\RB{\mathbb{R}}
\def\DB{\mathbb{D}}
\def\FB{\mathbb{F}}
\def\LB{\mathbb{L}}
\def\HB{\mathbb{H}}
\def\MB{\mathbb{M}}
\def\PB{\mathbb{P}}
\def\EB{\mathbb{E}}
\def\t{\tau}

\def\normep#1{\left|{#1}\right|_p}
\def\Tr#1{{\rm Tr}\left[#1\right]}

\def \Sum{\displaystyle\sum}
\def \Prod{\displaystyle\prod}
\def \Int{\displaystyle\int}
\def \Frac{\displaystyle\frac}
\def \Inf{\displaystyle\inf}
\def \Sup{\displaystyle\sup}
\def \Lim{\displaystyle\lim}
\def \Liminf{\displaystyle\liminf}
\def \Limsup{\displaystyle\limsup}
\def \Max{\displaystyle\max}
\def \Min{\displaystyle\min}
\newcommand{\dproof}{\noindent {Proof.} \quad}
\newcommand{\fproof}{\hfill $\square$ \bigskip}
\endlocaldefs
\begin{document}







\begin{frontmatter}
	\title{Non linear optimal stopping problem and Reflected BSDE in the predictable setting}
	\runtitle{Optimal Stopping in General Predictable Framework}
	
	\begin{aug}
		\author{\fnms{Siham} \snm{Bouhadou}\thanksref{m1}\ead[label=e1]{sihambouhadou@gmail.com}}
		\and
		\author{\fnms{Youssef} \snm{Ouknine}\thanksref{m1,m2}\ead[label=e2]{ouknine@uca.ac.ma} \ead[label=e3]{youssef.ouknine@um6p.ma}}

\runauthor{ Bouhadou and Ouknine}

		\affiliation{Cadi Ayyad University\thanksmark{m1} and Mohammed VI Polytechnic University\thanksmark{m2}}
		
		\address{Department of Mathematics,\\ 
			Faculty of Sciences Semlalia,\\
			 Cadi Ayyad University,\\
			  Marrakech, Morocco.\\
			\printead*{e1}}
		
		\address{Department of Mathematics,\\ 
			Faculty of Sciences Semlalia,\\
			Cadi Ayyad University, Marrakech.\\
			Mohammed VI Polytechnic University\\
			Benguerir, Morocco. \\		
			\printead*{e2}\\
			\printead*{e3}
		}
	\end{aug}
\begin{abstract}
	\begin{spacing}{1.2}
	 In the first part of this paper, we study RBSDEs in the case where the filtration is non quasi-left continuous and the lower obstacle is given by a predictable process.  We prove the existence and uniqueness by using some results of optimal stopping theory in the predictable setting, some tools from general theory of processes as the Mertens decomposition of predictable strong supermartingale. 
	In the second part we introduce an optimal stopping problem indexed by predictable stopping times with the non linear predictable $g$ expectation induced by an appropriate BSDE. We establish some  useful properties of ${\cal{E}}^{p,g}$-supremartingales.  Moreover,  we show the existence of an optimal predictable stopping time, and we characterize the predictable value function in terms of the first component of RBSDEs studied in the first part.  
\end{spacing}
\end{abstract}

\begin{keyword}
\kwd{optimal stopping}
\kwd{no quasi left continuous filtration}
\kwd{ predictable supermartingale}
\kwd{american options}
\kwd{predictable Snell envelope}
\kwd{predictable optimal stopping time}
\kwd{american options}
\kwd{reflected backward stochastic differential equation}
\end{keyword}



\end{frontmatter}
\section{Introduction}\label{sec1}
\begin{spacing}{1.2}
The notion of nonlinear Backward stochastic differential equation BSDE was introduced by Pardoux and Peng in the seminal work \cite{Pape90} when the noise is driven by a Brownian motion. A solution of this equation associated with a terminal value $\xi$ and a driver $g(t,\omega,y,z)$ is a couple of stochastic processes $(Y,Z)$ living in appropriate spaces, such that 
$$Y_t=\xi+ \int_t^T g(s,Y_s,Z_s)ds-\int_t^TZ_s dW_s$$
for all $t\leq T$. Where $W$ is a Brownian motion and the processes are adapted to natural filtration of $W$. In \cite{Pape92}, they provide Feynman Kac formulation representations of  solutions of non linear parabolic differential equations. Since then, these equations have found numerous applications in many fields of mathematics such as finance (see e.g. \cite{barr:elka:05,Rouge}), stochastic optimal control and games (see e.g. \cite{H-LP,H}), or partial differential equations (see e.g. \cite{Pape92}).

The theory of BSDEs has been extensively studied in the context of a filtration which is generated by Brownian motion, possibly with the addition of Poisson jumps (see e.g. \cite{HO1}, \cite{HO2}). In the case of more general filtration, one needs to introduce another component  in the above definition, namely a martingale $M$ that it is orthogonal to $W$:
	$$X_t=\xi+\int_{t}^T g(s, X_s,  Z_s)ds-\int_{t}^T  Z_s dW_s-M_{T}+M_{t} \text{   for all } t\in[0,T] \text{ a.s. }  $$
	These equations were first introduced by El Karoui and Huang \cite{elka:huan:9}, and treated recently by Kruse and Popier \cite{popier} to handle more general filtrations, which are not necessarily  quasi left continuous.   We remind the reader that the filtration $\mathbb{F}=(\mathcal{F}_t)$, assumed to satisfy the usual hypotheses, is said to be quasi left continuous if , for any predictable stopping time $\t$, one has ${\cal F}_{\t}={\cal F}_{\t^-}$. Intuitively, this means that martingales with respect to $\mathbb{F}$ cannot have predictable times of jumps. To understand the difficulty induced by avoiding the quasi-left-continuity assumption, we refer the reader to the work \cite{bouchard}. A huge part of the literature studied BSDEs
	 in the context of  quasi-left-continuous filtrations. Thus, the attention has been given to BSDEs when the stochastic terminal value $\xi$ is in $L^2({\cal{F}}_T)$, where $T$ is a fixed finite terminal time, for which the solution is required to be adapted to the natural filtration.
In the general setting, it seems natural to ask what's can be the formulation of BSDEs if we take $\xi$ is in $L^2({\cal{F}}_{T^-})$ and if we wish that the first component of BSDE  $X$ satisfy: $$X_\t\;\;\;\mbox{is}\;\;\;\;\;\;{\cal{F}}_{\t^-}\mbox{-measurable for each predictable stopping time}\;\;\; \t.$$ In this paper, we  show that the formulation can be as follows 
\begin{equation}\label{quas}
 X_t=\xi+\int_{t}^T g(s, X_s,  \pi_s)ds-\int_{t}^T  \pi_s dW_s-M_{T^-}+M_{t^-} \text{   for all } t\in[0,T] \text{ a.s. }  
\end{equation}
 Where $M_{t^-}$ denote the  left limit of the martingale $M$ at $t$. Note that Contrary to the classical case, these BSDEs are able to deal with any situations where martingales can jump at a predictable time with positive probability. 
 	A significant use of these  equations  is to generate  a  new family of "non linear expectations" or "nonlinear evaluations", which we  call \emph{predictable} $g$-conditional expectation  defined in the same spirit of \cite{Pe04}. In the present paper, we are  interested in a generalisation of classical optimal stopping problem where the linear expectation is replaced by the predictable $g$ expectation. 
 \\
 Optimal stopping problems with one agent whose payoff is assessed by a non-linear expectation has been introduced in El Karoui and Quenez \cite{EQ96} in the case of a Brownian filtration and a continuous pay-off process $\xi$. The problem has been generalized to the case of a Wiener-Poisson filtration and a right-continuous pay-off process $\xi$ in Quenez and Sulem \cite{QuenSul2}. The paper of  Grigorova et al. \cite{MG} is the first to consider the case of  non-right-continuous  pay-off  process $\xi$  by assuming  the weaker  assumption of right- uppersemicontinuity.  In \cite{Imkeller2},  Grigorova et al. study  this optimal stopping problem \ without making any regularity assumptions on process $\xi$.\\
 From practical point of view, the agent would like to choose his strategy in such way to have a minimal risk as possible. In our setting we will show that it is possible to fulfill these desire, by studying the risk measure induced by the BSDEs \eqref{quas} with a Lipschitz driver $g$ in the predictable setting.
The optimal stopping in this context can be formulated as follows: given a dynamic financial position process $\xi$, represented by a ladlag predictable process, we want to determine a predictable stopping time which minimizes the risk of position $\xi$. For a predictable stopping time $S$ such that  $0\leq S\leq T$ a.s. (where $T>0$ is a
fixed terminal horizon) , we define
\begin{equation}\label{first} 
V_p(S): =   {\rm ess} \sup_{\t \in \stops}{\cal E}^{p,g}_{S,\t}(\xi_{\t}),
\end{equation}
 where  $\mathcal{T}_{S}^p$
 denotes the set of predictable stopping times valued a.s. in $[S,T]$ and ${\cal E}^{p,g}_{S, \t}(\cdot)$ denotes  the predictable $g$-conditional expectation.  
 \\
 The study of classical optimal stopping problem in the predictable framework (corresponding to $g=0$ in  \eqref{first}), dates back to El Karoui  in the work \cite{EK}, in which she mentioned the complexity to exhibit conditions ensuring the existence of a solution by using the penalization method of Maingueneau \cite{Maingueneau}. In this work, we focus on problem \eqref{first}  where, on one hand, the filtration
is not quasi-left continuous, on the other
hand, the reward process $\xi$ is not right continuous but only a ladlag predictable process, by using its links with an appropriate RBSDE.
 \\
Let us recall that RBSDEs have been introduced  by El Karoui et al. \cite{ElKaroui97}  and have proved useful, for instance, in the study of American options. The work by El Karoui et al. \cite{ElKaroui97} considers the case of a Brownian filtration and a continuous obstacle.   There have been several extensions of this work to the case of a discontinuous obstacle (cf. \cite{Ham}, \cite{HO1}, \cite{Essaky}, \cite{HO2}, \cite{QuenSul2}), or to the no right continuous  case in \cite{MG}, or to more general framework, without any regularity assumption in \cite{Imkeller2}.
In the first part of the present paper, we formulated RBSDEs in predictable setting, where the obstacle is assumed to be a completely irregular predictable process. We prove here the existence and uniqueness of the solution  by using some results of optimal stopping Problem of \cite{Karoui}, some tools from the general theory of processes .  More precisely, we use  Merten's decomposition of predictable strong supermartingale (cf. Meyer \cite{Meyer_cours}), which can be seen as a generalization of Doob-Meyer decomposition, and a suitable version of Gal'chouk-Lenglart's formula (see \cite{Galchouk}).
In the second part of the paper, we  begin by
 providing some additional properties of the operator ${\cal{E}}^{p,g}$, which are specific to the predictable framework. We also prove that the solution of RBSDEs studied in the first part, is the value of the non linear problem \eqref{first}.
 We show  also the existence of a predictable stopping time for the problem  under some additional conditions on $\xi$.\\
Now, Let us outline some differences in our paper compared to existing literature.  First, we restrict our attention to that
  the solution of BSDE \eqref{quas}  is no longer right continuous.  Moreover, the non linear optimal stopping problem considered in our setting is formulated only over predictable stopping times.\\
 Compared to \cite{MG}, \cite{Imkeller2},  the process of right limits of the additional increasing process $A+B_{-}$, which pushes the first component of the solution of  RBSDE,  to stay above the predictable obstacle $\xi$ is predictable. Moreover, the role of $B$ is to make necessary jumps to keep the solution $Y$ to stay above the barrier and it doesn't  act only when $Y$ has right jumps which is the case in \cite{MG} and \cite{Imkeller2}, but it acts also at predictable times $\t $ for which $Y_\t \neq {}^pY_{\t}^+$.\\
\textbf{Organisation of the paper:}
\\
 The paper is organized as follows: the second section is dedicated to some preliminary definitions and properties. In section \ref{section-RBSDE}, we define reflected  BSDEs in the predictable framework, we prove also the existence and uniqueness of the solution. In Section \ref{section-optimal}, we derive some useful properties of the operator ${\cal E}^{p,g}$. In Section \ref{Existence of predictable}, we prove the existence of 	a predictable optimal stopping time under an additional assumption on the obstacle $\xi$. We characterize the value function in terms of the first component of the solution of  RBSDEs studied in the first part.  In Section \ref{section snell}, we prove some additional results on the strong predictable Snell envelope. In the appendix, we recall some tools from the general theory of processes  as Merten's decomposition for predictable supermartingale  and Galchuk-Lenglart's formula.

 \section{Preliminaries and definitions}
 We start with some notations. We fix a stochastic base with finite horizon $T \in \mathbb{R}_{+}^*$, $(\Omega,\mathbb{F}=(\mathcal{F}_{t})_{t \in {[0,T]}},P).$ 
 We assume that the filtration $\mathbb{F}$ satisfies the usual assumptions of right continuity and completeness Importantly, we  assume that the filtration is not quasi-left continuous. Let $W$ is one-dimensional $\mathbb{F}$-Brownian motion. \\
  We recall that a stopping time $\t$ is called predictable if there exist a sequence $(\t_n \leq \t)$ that are strictly smaller than $\t$ on$\{\t>0\}$ and increase to $\t$ a.s.
 		\begin{itemize}
 			\item We denote by $\stopo$ the collection of all \emph{predictable} stopping times  $\tau$ with values in $[0,T]$.
 			More generally , we denote $\stops$ (resp. $\mathcal{T}_{S^+}^p$) the class of predictable  stopping times $\tau\in \stopo$
 			with $S\leq\tau$ a.s.  (resp. $\tau>S $ a.s. on $\{S<T\}$ and $\tau=T$ a.s. on $\{S=T\}$).
 		\end{itemize}
	We use also the following notation:
	\begin{itemize}
		\item ${\cal P}$ (resp. $\mathcal{O}$) is  the predictable (resp. optional) $\sigma$-algebra
		on $ \Omega\times [0,T]$.
		\item $\mathrm{Prog}$ is the progressive $\sigma$-algebra
		on $ \Omega\times [0,T]$.
		\item ${\cal B}(\R)$ (resp. ${\cal B}(\R^2)$)   is the Borel
		$\sigma$-algebra on $\R$ (resp. $\R^2$).
		\item
		$L^2({\cal F}_{T^-})$  is the set of random variables which are  $\Fc
		_{T^-}$-measurable and square-integrable.

		\item    $\H^{2}$ is the set of
			$\R$-valued predictable processes $\xi$ with
				$\|\xi\|^2_{\H^{2}} := E \left[\int_0 ^T |\xi_t|^2dt\right]< \infty$

			\end{itemize}
		


		We denote by  $ {\cal S}^{2,p}$ the vector space of
		$\R$-valued predictable (not necessarily cadlag)
		processes $\xi$ such that
		$\vertiii{\xi}^2_{{\cal S}^{2,p}} := E[\esssup_{\tau\in\stopo} |\xi_\tau |^2] <  \infty.$ 
		By using similar arguments as in the proof of  Proposition 2.1 in \cite{MG},  one can show that the mapping $\vertiii{\cdot}_{{\cal S}^{2,p}}$ is a norm on the space $ {\cal S}^{2,p}$,  and ${\cal S}^{2,p}$ endowed with this norm is a Banach space.\\
		$\bullet$ Let ${\cal M}^{2}$  be the set of square integrable  martingales $M= (M_t)_{t\in[0,T]}$ with $M_0=0$. 
\\
		$\bullet$ Let ${\cal M}^{2,\bot}$ be the subspace of martingales $M \in {\cal M}^{2}$
		satisfying  $\langle M, W \rangle_\cdot =0$.

		%

For a ladlag process $X$, we denote by $X_{t+}$ and $X_{t-}$ the right-hand and left-hand limit of $X$ at  $t$. We denote by $\Delta_+ X_t:=X_{t_+}-X_t$ the size of the right jump of $X$ at $t$, and by   $\Delta X_t:=X_t-X_{t-}$ the size of the left jump of $X$ at $t$.
Let us recall  the key section theorem related to indistinguishability of optional processes or predictable processes.

\begin{theorem}[Section Theorem]\label{section}
Let $X=(X_t)$ and  $Y=(Y_t)$ be two optional (resp. predictable) processes. If for every bounded stopping time (resp. predictable time) $\t$, we have $X_\t\leq Y_\t$ a.s. (resp. $X_\t=Y_\t$ a.s.), then $X\leq Y$ (resp. X and Y are indistinguishable).
\end{theorem}
Let us recall the following orthogonal decomposition property of martingales in ${\cal M}^2$.
		\begin{lemma}[Lemma III.4.24 in \cite{JS}]\label{lemma1}
			For each $N \in {\cal M}^2$, there exists a unique couple $(Z,M) \in \HB^2 \times  {\cal M}^{2,\perp} $
			\begin{equation}
			N_t =\int_0^t Z_s dW_s+M_t. \; \; \forall t \in [0,T]\;\;\;\mbox{a.s.}
			\end{equation}
		\end{lemma}
\begin{definition}[Driver, Lipschitz driver]\label{defd}
			A function $g$ is said to be a {\em driver} if
			\begin{itemize}
				\item (measurability)
				$g: \Omega  \times [0,T]  \times \R^2  \rightarrow \R $\\
				$(\omega, t,y)\mapsto g(\omega, t,y, z)  $
				is $ \mathcal{P} \otimes {\cal B}(\R^2) 
				- $ measurable,
				\item (integrability) $E[\int_0^T g(t,0,0)^2dt] < + \infty$.
			\end{itemize}
			A driver $g$ is called a {\em Lipschitz driver} if moreover there
			exists a constant $ K \geq 0$ such that $dP \otimes dt$-a.e.\,, for
			each $(y_1, z_1)\in \R^2$, $(y_2, z_2)\in \R^2$,
			$$|g(\omega, t, y_1, z_1) - g(\omega, t, y_2, z_2)| \leq
			K (|y_1 - y_2| + |z_1 - z_2|).$$
			A pair $(g,\xi)$ such that $g$ is a Lipschitz driver and $\xi \in L^2(\cf_{T}) $   is called a pair of
			standard data, or a pair of standard parameters.
		\end{definition}
		\begin{definition}[BSDE, conditional  predictable $g$-expectation]\label{BSDE}
		Let  $g$ be a Lipschitz driver, and $\xi $ in $L^2(\cf_{T^-})$. We will formulate the  BSDE associated with Lipschitz driver $g$, terminal time $T$, and terminal condition $\xi$,  as follows:
		$$X_t=\xi+\int_{t}^T g(s, X_s,  \pi_s)ds-\int_{t}^T  \pi_s dW_s-M_{T^-}+M_{t^-}\;\; \text{   for all } t\in[0,T] \text{ a.s. }  $$
		We recall that the above BSDE in the case of quasi-left continuous filtration and when the solution is required to be just adapted, corresponds to the  classical BSDE with standard data which has been widely considered in the literature. In this paper, we will show that the above BSDE admits a unique solution  $(X,\pi,M)\in{\cal S}^{2,p} \times \H^2 \times  {\cal M}^{2,\perp} $ for which the first component $X$ is left continuous.\\ For $t\in[0,T]$, we introduce  the (non-linear) operator $\mathcal{E}^{p,g}_{t,T}(\cdot): L^2(\cf_{T^-})\rightarrow L^2(\cf_{t^-})$ which maps a given terminal condition
		$\xi\in L^2(\cf_{T^-})$ to the position $X_t$ (at time $t$) of the first component  of the solution of the above BSDE 
	, and we will call it \emph{ predictable conditional $g$-expectation at time $t$}.  As usual, this notion can can be extended to
		the case where the (deterministic) terminal time 
		$T$
		is replaced by a (more general) predictable stopping time $\tau\in\stopo$,  $t$ is replaced by a predictable stopping time $S$ such that $S\leq \tau$ a.s. and  the domain $L^2(\cf_{T^-})$ of the operator    is replaced  by $L^2(\cf_{\tau^-})$.
		
	\end{definition}
	\section{Reflected BSDE  whose obstacle is predictable in the case of non quasi-left continous filtration}\label{section-RBSDE}
	The notion of reflected BSDEs has been recently  studied by Grigorova, Imkeller, Ouknine and Quenez in the  seminal paper \cite{Imkeller2}, in the case of an optional completely irregular obstacle and a general filtration.  We began this section, by giving a formulation of this notion of reflected BSDEs  in the predictable setting, for which the solution of such equations is predictable and constrained to be greater than a given irregular predictable process called \emph{predictable irregular obstacle}. \\
	Let $T$ be a fixed terminal time and $g$ be a Lipschitz driver. Let $\xi$ be a ladlag predictable process in ${\cal S}^{2,p}$, called obstacle or barrier in ${\cal S}^{2,p}$.
\begin{definition}\label{defRBSDE} 
	A process $(Y,Z,M,A,B)$ is said to be a predictable solution to the reflected BSDE with parameters $(g,\xi)$, where $g$ is a driver and $\xi$ is a predictable obstacle, if 
\begin{eqnarray}\label{prince}
\left\{
\begin{array}{ll}
(i)&(Y,Z,M,A,B)\in {\cal S}^{p,2} \times \H^{2} \times {\cal M}^{2,\perp} \times {\cal S}^{2,p}\times {\cal S}^{2,p} \;\;\text{and a.s. for all} \;\;\tau\in\stopo
\nonumber\\

(ii)&Y_\tau=\xi_T+\int_{\tau}^T g(t,Y_t,  Z_t)dt-\int_{\tau}^T  Z_t dW_t-M_{T^-}+M_{\t^-} +A_T-A_\tau+B_{T-} -B_{\tau-} \text{ a.s.}, \\ \label{RBSDE}
(iii)&  Y_\t \geq \xi_\t  \text{ a.s.,} \text{ for all } \t\in\stopo,\\ \label{RBSDE_inequality_barrier}\nonumber
(iv)&  A \text{ is a nondecreasing right-continuous predictable process with} \;A_0= 0 \text{ and such that } \\
&\int_0^T {\bf 1}_{\{Y_t > \xi_t\}} dA^c_t = 0 \text{ a.s. and }\;(Y_{\tau-}-\xi_{\tau-})(A^d_{\tau}-A^d_{\tau-})=0 \text{ a.s. for all (predictable) }\tau\in\stopo, \label{RBSDE_A}\\   
(v)& B \text{ is a nondecreasing right-continuous \textbf{ predictable} purely discontinuous} \\
&\text{ process with } B_{0-}= 0, \text{ and such that }

(Y_{\tau}-\xi_{\tau})(B_{\tau}-B_{\tau-})=0 \text{ a.s. for all }\tau\in\stopo. \label{Min-C}
\label{RBSDE_C}
\end{array}
\right.
\end{eqnarray}
\end{definition}
Here  $A^c$ denotes the continuous part of the non-decreasing process $A$ and $A^d$  its discontinuous part.
\begin{remark}
	The equations (iv) and (v) are referred to as minimality conditions or Skorohod conditions.
\end{remark}
\begin{remark}
	In Definition \ref{defRBSDE}, we have given the notion of reflected BSDE with lower reflecting predictable barrier. However one could given the notion of reflected BSDE with upper predictable reflecting barrier.\\
	A quintuple $(Y,Z, M,A,B)$ is s solution for a reflected BSDE with  predictable lower reflecting  barrier $\xi$ and  a driver $g$  iff $(-Y,-Z, -M,A,B)$ is a solution for the reflected BSDE with a predictable upper reflecting  barrier associated with $-g(t,w,-y,-z)$.
\end{remark}
\begin{remark}
	The role of $A^d$ is to make necessary jumps to keep $Y$ above the barrier, and it always acts when $Y$ has left jumps, which occurs at left jumps of $\xi$.
	The process $A^c$ does act only when the process $Y$ reaches $\xi$ either at its continuity.
\end{remark}
\begin{remark}\label{jumpsA}
	If we rewrite the equation (ii)  forwardly, we can see that the left jumps of the process $Y$  verify $Y_\t-Y_{\t^-}=-(A_\t-A_{\t^-})$ a.s. for each predictable stopping time $\t \in \stopo$. Moreover, the jump processes  satisfy $\Delta Y\equiv\Delta A$. Indeed, the processes $ Y$ and $A$ are predictable, thus $\Delta Y$ and $\Delta A$ are also predictable. The result follows from an application of section theorem (see Theorem \ref{section}).
\end{remark}
\begin{remark}
	In our framework the filtration is not quasi-left continuous, the martingales have  totally inaccessible jumps and can also jump at predictable times.	
\end{remark}
\begin{remark}
	We restrict our attention to the fact that The term  $M_-$ in the equation satisfied by $Y$ is not a martingale but the predictable projection of the martingale $M$ (see Corollary \ref{cor-predictable} in the Appendix) . 
\end{remark}
	\begin{remark}\label{jumps}		
If $(Y,Z,M,A,B)$ is a solution of RBSDE defined above, and if $\t \in \stopo$, then $\Delta B_{\t}=Y_\t-{}^pY^+_{\t}$. This equality follows from the equation (ii) and the fact that $M$ is right continuous, ${}^pM=M_-$  and the fact that both $A$ and $B$ are predictable processes.
Moreover, the  processes $(Y_t-{}^pY^+_{t})_{t \in [0,T]}$ and $(\Delta B)_{t \in [0,T]}$ are indistinguishable. This is due to the fact that   $ Y$, ${}^pY^+$ and $B$ are all predictable.
\end{remark}
\begin{remark}\label{rem upper}
	Note also that $Y\geq {}^pY^+$ up to an evanescent set. 
The proof of this claim is due to the fact of non-decreasingness of (almost all trajectories of) $B$ and similar arguments as in Remark \ref{jumps}.
\end{remark}
Now, let us recall the definition of the predictable strong supermartingale. This notion has introduced by Chung and Glover \cite{Bhung}, see Appendix I of the book of Dellacherie and Meyer \cite{DM2} for the subsequent
concept.
\begin{definition}Let $(Y)_{t \in [0,T]}$ be a real valued process. We say that $Y$ is a  predictable strong supermartingale process if 
	\begin{itemize} 
		\item  $Y$ is predictable;
		\item $Y_\tau $ is integrable for all $\tau \in \mathcal{T}_0^p$.
		\item for all predictable stopping times $S \leq  \t$ $$Y_S \geq E[Y_\tau\mid \Fc_{S^-}].$$
	\end{itemize}
	
\end{definition}
\begin{remark}\label{rem.ladlag}
	Every predictable strong supermartingale is indistinguishable from a ladlag process,  see \cite{DM2}.
\end{remark}
\begin{remark}\label{strong}		
	If $(Y,Z,M,A,B)$ is a solution of RBSDE defined above, then the processes $(Y_t+\int_0^t g(Y_s,Z_s)ds)_{t \in [0,T] }$ and  $(^{p}Y_{t}^++\int_0^t g(Y_s,Z_s)ds)_{t \in [0,T] }$ are predictable strong supermartingale processes.
\end{remark}

Now, we state a priori estimates on the solutions.
\begin{lemma}\label{estimation}
	Let $(Y^1,Z^1,M^1,A^1,B^1) \in ({\cal S}^{2,p}\times \HB^2 \times {\cal M}^{2,\perp} \times {\cal S}^{2,p}\times {\cal S}^{2,p})$ (resp.  $(Y^2,Z^2,M^2,A^2,B^2) \in ({\cal S}^{2,p}\times \HB^2 \times {\cal M}^{2,\perp} \times {\cal S}^{2,p}\times {\cal S}^{2,p})$ )  be a solution to the RBSDE associated with data $(g^1,\xi^1)$  (resp., $(g^2,\xi^2)$). We set 
	$\tilde Y:=Y^1-Y^2$, $\tilde{g}(\omega,t):=g^1(\omega,t)-g^2(\omega,t)$.\\
	There exists $c>0$ such that for all $\varepsilon>0$,  for all $\beta\geq\frac 1 {\varepsilon^2}$ we have 
	\begin{align}\label{eq_initial_Lemma_estimate}
	\|\tilde Z\|^2_\beta\leq \varepsilon^2\|\tilde g\|^2_\beta, \,\,  \text{  and  } \vvertiii{\tilde Y}^2_\beta \leq 2\varepsilon^2(1+8c^2) \|\tilde g\|^2_\beta.
	\end{align}
\end{lemma}
\dproof The proof is given in the Appendix.
\fproof
\begin{remark}
	The  estimates given in Lemma \ref{estimation} still hold in case of predictable non reflected BSDE in general filtration  (see Definition \ref{BSDE}).
\end{remark}

We  introduce the following definition.
\begin{definition} \label{defr} A progressive process $(\xi_t)$ (resp. integrable) is said to be {\em left-upper semicontinuous (l.u.s.c.) along stopping times} (resp. along stopping times in expectation ) if for all $\tau \in {\cal T}_0$ and for each non decreasing sequence of stopping times $ (\tau_n)$ such that $\tau^n \uparrow \tau$ a.s.\,,
		\begin{equation}\label{usc}
		\xi_{\tau} \geq \limsup_{n\to \infty} \xi_{\tau_n} \quad \mbox{a.s.} \quad \text{ (resp. } E[\xi_{\tau}] \geq \limsup_{n\to \infty} E[\xi_{\tau_n}] ).
		\end{equation}
	\end{definition}
	\begin{remark}\label{LC}
		Note that  when $(\xi_t)$ is left-limited, then $(\xi_t)$ is {\em left-upper semicontinuous (l.u.s.c.) along stopping times} if and only if for all predictable stopping time $\tau \in \stopo$, 
		$\xi_{\tau} \geq  \xi_{\tau_-}\quad\mbox{a.s.}$ In particular, if  $(\xi_t)$ is left-limited predictable process, $(\xi_t)$ is l.u.s.c if and only if  $\xi \geq \xi_{-}$ up to an evanescent set.
	\end{remark}

		\begin{remark}\label{remarkA}
		As a direct consequence of the equation satisfied by $Y$, we have  $Y_\t-Y_{\t^-}=-(A_\t-A_{\t^-})$ a.s. for each predictable stopping time $\t \in \stopo$. Thus, $Y_{\t^-} \geq Y_{\t}$ a.s. for all  $\t \in \stopo$. By the section theorem $Y_{-} \geq Y$ up to an evanescent set. 
	\end{remark}
	\begin{lemma}\label{lem left}
		If $\xi$ is predictable l.u.s.c  along stopping times, then $Y$ is left continuous. On other words, the process $A$ is continuous.
			\end{lemma}
\dproof
		The  Remark \ref{remarkA} combined with the fact that $Y \geq \xi$ and that $\xi$ is l.u.s.c along stopping times lead to $Y_{\t^-} \geq Y_{\t}\geq   \xi_{\t} \geq \xi_{\t^-} $.  Since,  $Y_{\t^-}\geq \xi_{\t^-}$, we have two cases. If $\t$ is such that $Y_{\t^-}=\xi_{\t^-}$, then the above inequalities become equalities and we obtain  $Y_{\t^-}= Y_{\t}$. If
		$\t$ is such that $Y_{\t^-}>\xi_{\t^-}$, then $Y_{\t}-Y_{\t^-}=-(A_{\t}-A_{\t^-})=0$ [due to Remark \ref{jumpsA} and to
		minimality condition (iv) satisfied by $A$]. Thus, in both cases, $Y_{\t}=Y_{\t^-}$, which proves the left continuity of $Y$.
	\fproof
		%
\begin{remark}\label{rem-jumps}
If $Y$ is right continuous, then $B$ is indistinguishable from the null process $0$. In fact, $Y$ is predictable. Moreover, it is right continuous by hypothesis. Thus, $^pY^+_{\t}=Y_\t$ for all $\t \in \stopo$. We get by Remark \ref{jumps} that the jump process $\Delta B$ is null.  we derive the desired result from the fact that $B$ is a non-decreasing right-continuous
adapted purely discontinuous process with  $B_{0-}= 0$.				
\end{remark}
		\begin{lemma}\label{lem right}
		If $\xi$ is  right-continuous  predictable obstacle and if the filtration is quasi-left continuous, then $Y$ is right-continuous.
		\end{lemma}
	\dproof
		Indeed, we have by Remark \ref{rem upper}, $Y_t \geq {}^p Y^+_{t}\geq  {}^p\xi_{t}^+ $. Since $\xi$ is right continuous predictable process, $ {}^p\xi_{t}^+=\xi_{t}$.  Hence,  $Y_t \geq {}^p Y^+_{t}\geq\xi_{t} $. If $t$ is such that $Y_t=\xi_t$, then $Y_t={}^p Y^+_{t}
		$. If
		$t$ is such that $Y_t > \xi_t$, then $Y_t-{}^pY^+_{t} =B_t-B_{t^-} = 0$ [due to Remark \ref{jumps} and to
		minimality condition (v) satisfied by $B$]. Thus, in both cases, $Y_t = {}^p Y^+_{t}$. Since the filtration is quasi-left continuous, the martingale $M$ is quasi-left continuous. By the same arguments as in the Remark \ref{rem-jumps},  the process $B$ is indistinguishable from the null process $0$. Therefore, Y is right-continuous.
\fproof	
\begin{remark}
	If $\xi$ is  right-continuous l.u.s.c  predictable process and if the filtration is quasi-left continuous,  then, $Y$ is continuous. This is a direct consequence of Lemmas \ref{lem left} and \ref{lem right} .
\end{remark}

\subsection{\textbf{Existence and uniqueness of the solution}}

We are now going to prove That RBSDE from Definition \ref{defRBSDE} has a unique solution, by using Snell's method envelope. In a first step, we suppose that $g$ does not depend on $(y,z)$.
	\begin{theorem}\label{existence-linear}
		Consider a couple $(g,\xi)$ where $g$
 is a progressive process and suppose that $g(\omega,t,y,z)=g(\omega,t)$ does not  depend on $(y,z)$ with  $E[\int_0^T g(t)^2dt]< \infty$. Then, there exists a unique solution $(Y,Z,M,A,B) \in {\cal S}^{2,p} \times  \HB^2 \times {\cal M}^{2,\perp}\times {\cal S}^{2,p}\times {\cal S}^{2,p}$ of the RBSDE from Definition \ref{defRBSDE}  , and for each $S \in \stopo$, we have:
 \begin{eqnarray*}
 	Y_S= \esssup_{\t \in \stops} E\left[ \xi_\tau + \int_S^\t g_u du \mid \Fc_{S^-}\right].
 \end{eqnarray*}
Moreover, the following  properties hold:
\begin{description}
	\item[(ii)]	We have	$ Y_S=\xi_S\vee {}^pY^+_{S}$ a.s.
	for all $S \in \stopo$.
	\item[(iii)]	We have	$ Y_{S^-}=\xi_{S^-}\vee Y_{S}$ a.s	for all $S \in \stopo$.
	\item[(iv)] For each $S \in  {\cal T}_{0}^p$ and for each $\alpha \in ]0,1[$, we set
	\[\tau^{\alpha}(S):= \inf\left\{t\geq S \,, \alpha Y_t+(\alpha-1)\int_0^tg(u)du\leq \xi_t \right\}.\]
	Then, 
	for each $\alpha  \in ]0,1[$ and for each $S \in \stopo$, 
	\begin{equation}\label{eq.mar}
	Y_S=E\left[-\int_0^{S} g(u)du+\int_0^{\tau^\alpha(S)} Z_s dW_s+M_{\tau^\alpha(S)}-A_{\tau^\alpha(S)}-B_{\tau^\alpha(S)^{-}}|\mathcal{F}_{S^-}\right]\;\;\;\; \mbox{a.s.}
	\end{equation}
	Moreover,
	$A_{ \tau^\alpha(S)}=A_S$ and  $B_{\tau^\alpha(S)^-}=B_{S^-}$.	
\end{description}

 \end{theorem}
The proof of the Theorem \ref{existence-linear} relies on some lemmas.\\
First, let introduce 
 the value function $Y(S)$  defined at each $S \in \stopo$ by
\begin{eqnarray}\label{ddeux-2}
Y  (S):= \esssup_{\t \in \stops} E\left[ \xi_\tau + \int_S^\t g_u du \mid \Fc_{S^-}\right].
\end{eqnarray}
For the reader's convenience, we give a detailed study of  the predictable value function in section \ref*{section snell}.
\begin{lemma}\label{aggreg}	
	 There exists a ladlag predictable process $( Y_t)_{t\in[0,T]}$ which aggregates the family
		$(Y(S))_{S\in\stopo}$ (i.e. $ Y_S= Y (S)$ a.s. for all $S\in\stopo$).\\
		Moreover, the process $(U_t:= Y_t + \int_0^t g_u du)_{t\in[0,T]}$ is characterized as the predictable Snell envelope associated with the process 
		$( \xi_t + \int_0^{t}g_u du)_{t \in [0,T]}$, that is the
		smallest predictable supermartingale
		greater than or equal to the process  $( \xi_t + \int_0^{t}g_u du)_{t \in [0,T]}$.
		\end{lemma}
	\dproof
	Let us prove the first assertion. For each $S \in \stopo$, we  define the random variable $U(S)$ by 
	$$U(S):=Y(S)+\int_0^S g(u)du=\esssup_{\t \in \stops} E\left[ \xi_\tau + \int_0^\t g_udu \mid \Fc_{S^-}\right].$$
	Put $\psi_\t=\xi_\tau + \int_0^\t g(u)du$. 
	Since the process $(\xi_t)_{t\in [0,T]}$ is predictable, we get also that the process  $(\psi_t)_{t\in [0,T]}$ is predictable. Hence, the family $(U(S))_{S\in\stopo}$ is the predictable Snell envelope system associated to $(\psi_t)_{t\in [0,T]}$ (see second assertion of Lemma  \ref{lemm4}). Moreover, by a result of El Karoui \cite{Karoui} (see (i) Lemma  \ref{pro}  in the Appendix), there exists a strong predictable supermartingale process (which we denote again by $U$) which is the predictable Snell envelope associated to $\psi$ such that $U_S=U(S)\, \rm{ a.s.}$ for all $S\in\stopo$.  Thus, we have $$Y(S)=U(S)-\int_0^S g_u du=U_S-\int_0^S g_u du\, \;\;\rm{ a.s.} \;\; \mbox{for all}\; S\in\stopo .$$  On the other hand, by  Remark \ref{rem.ladlag}, almost all trajectories of $U$ are ladlag. Thus, we get that the ladlag predictable process $(Y_t)_{t\in[0,T]}=(U_t-\int_0^t g_u du)_{t\in[0,T]}$ aggregates the family $(Y(S))_{S\in\stopo}$. This yields  the desired result.
	\fproof
	\begin{lemma}
		The process $Y$ defined above verify the following statements:
				\begin{description}
\item[(i)]	We have	$ Y_S=\xi_S\vee {}^pY^+_{S}$ a.s.
			for all $S \in \stopo$.
\item[(ii)]	We have	$ Y_{S^-}=\xi_{S^-}\vee Y_{S}$ a.s. for all $S \in \stopo$.
		\end{description}
	\end{lemma}
	\dproof
By Lemma \ref{aggreg}, $(Y_t + \int_0^t g_u du)_{t\in[0,T]}$ is the predictable Snell envelope associated with the process 
$( \xi_t + \int_0^{t}g_u du)_{t \in [0,T]}$. Thus,  by Lemma \ref{pro} (ii) in the Appendix, we obtain, for all $S \in \stopo$
\begin{eqnarray*}
	Y_S+\int_0^{S}g_u du&=&\left(\xi_S+\int_0^S g(u)du\right)\vee {}^p \left(Y_{.}^+ +\int_0^{.}g_u du\right)_S\\
&=&\left(\xi_S+\int_0^S g(u)du\right)\vee \left({}^p Y_{S}^+ +\int_0^{S}g_u du\right).
\end{eqnarray*}
This yields $(i)$.\\
Since the processes $\xi$ and $Y$ are left limited and by the continuity of the process $\left(\int_0^{t}g_u du\right)_{t \in [0,T]}$, we obtain by Lemma \ref{pro} (iii),
\begin{eqnarray*}
	Y_{S^-}+\int_0^{S}g_u du&=&\left(\xi_{S^-}+\int_0^S g(u)du\right)\vee  \left(Y_{S} +\int_0^{S}g_u du\right).\\
\end{eqnarray*}
	Hence, we get $(ii)$.
\fproof
	\begin{lemma}\label{Mertensvalue}
		\begin{description}
			\item[(i)]  The predictable process $( Y_t)_{t\in[0,T]}$ is in ${\cal S}^{2,p}$ and admits the following predictable Mertens decomposition:
			\begin{equation}\label{eqmert}
			Y_{\tau}=Y_{0}-\int_0^\t g(u)du+N_{\tau^-}-A_{\tau}-B_{\tau^-}\;\;\; \mbox{for all}\;\tau\in\stopo.
			\end{equation}
			where $N \in {\cal M}^2$ , $A$ is a non-decreasing
			right-continuous predictable process such that $ A_0= 0$ and	$E(A_T^2)<\infty$, and $B$ is a non-decreasing right-continuous
			predictable purely discontinuous process such that  $B_{0-}= 0$ and	$E(B_T^2)<\infty$.
			\item[(ii)]
			For each $\tau \in\stopo$, we have
			$\Delta B_{\t}= {\bf 1}_{\{ Y_\t= \xi_\t\}}\Delta B_{\t}$ a.s.\,
			\item[(iii)]  For each  $\t \in\stopo$, we have
			$\Delta A_{\t}= {\bf 1}_{\{{ Y}_{\t-}= \, \xi_{\t^-}\}}\Delta
			A_{\t}$ a.s.
			\item[(iv)] For each $S \in  {\cal T}_{0}^p$ and for each $\alpha \in ]0,1[$, we set
			\[\tau^{\alpha}(S):= \inf\left\{t\geq S \,, \alpha { Y}_t+( \alpha-1)\int_0^tg(u)du\leq \xi_t\right \}.\]
			then 
			for each $\alpha  \in ]0,1[$ and for each $S \in \stopo$, 
			\begin{equation}\label{eq.mar}
			Y_S=E\left[-\int_0^{S} g(u)du+N_{\tau^\alpha(S)}-A_{\tau^\alpha(S)}-B_{\tau^\alpha(S)^{-}}|\mathcal{F}_{S^-}\right]\;\;\;\; \mbox{a.s.}
			\end{equation}
			Moreover,
			$A_{ \tau^\alpha(S)}=A_S$ and  $B_{\tau^\alpha(S)^-}=B_{S^-}$	
		\end{description}	
	\end{lemma}
\dproof
Note that $( \xi_t + \int_0^{t}g_u du)_{t \in [0,T]}$ is predictable process in ${\cal S}^{2,p}$. This is due to the fact that $\xi \in {\cal S}^{2,p} $ and $g \in \H^2$.
By Lemma \ref{aggreg}, $(Y_t + \int_0^t g_u du)_{t\in[0,T]}$ is the predictable Snell envelope associated with the process 
$( \xi_t + \int_0^{t}g_u du)_{t \in [0,T]}$. Thus, the Lemma \ref{Mertensvalue}  corresponds to Lemma \ref{principal} in the Appendix.
\\
\fproof
\\
In the following, we give the minimality property concerning $A^c$ in our predictable setting and for a no quasi left continuous filtration. In the Brownian  and right continuous framework, the proof can be found in  Karatzas and Shreve \cite{KS2} and it is based on some analytic arguments, which has proven  to be an efficient tool to generalize this result to the right upper semicontinuous case in \cite{MG}, and recently  to a completely irregular $\xi$ and in general filtration framework in \cite{Imkeller2}.  This last result can be seen as a generalization of the minimality property in the Brownian setting and for an optional $\xi$ in \cite{nouveau}, where the authors used another type of arguments.
\begin{lemma}\label{lemc}
	The continuous part $A^c$ of $A$ satisfies the equality $\int_0^T 1_{\{Y_{t^-}>\xi_{t^-}\}}dA_t^c=0$ a.s.
\end{lemma}
\dproof
The proof is based on Lemma \ref{Mertensvalue} (iv) which yields that 	$A_{ \tau^\alpha(S)}=A_S$ for each $S \in \stopo$ and for each $\alpha \in ]0,1[$, and the same arguments used in the proof of Theorem D $13$ in \cite{KS2}. For the convenience, we refer also  to the proof of Lemma $3.3$ \cite{Imkeller2}.

\dproof[Proof of Theorem \ref{existence-linear} ]
	\begin{description}
	\item[(1)]\textbf{Terminal condition:}
		We get by aggregation equality (see Lemma \ref{aggreg},	) combined with equation \eqref{ddeux-2} that $Y_T=Y(T)=\xi_T$ a.s.
	\item[(2)] $\mathbf{(Y,Z,M, A,B)\in{\cal S}^{2,p}\times \HB^2 \times  {\cal M}^{2,\perp}\times {\cal S}^{2,p} \times {\cal S}^{2,p}} $ \textbf{and the first component $Y$ satisfies the equation (ii)   of Definition \ref{defRBSDE}}. In fact, By Lemma \ref{Mertensvalue}\;$ (i), ( Y_t)_{t\in[0,T]}$ is in ${\cal S}^{2,p}$ and admits the following predictable Mertens decomposition:
	\begin{equation}\label{eqmert}
	Y_{\tau}=Y_{0}-\int_0^\t g(t)dt+N_{\tau^-}-A_{\tau}-B_{\tau^-}\;\;\; \mbox{for all}\;\tau\in\stopo.
	\end{equation}
Moreover, $(A,B) \in {\cal S}^{2,p} \times {\cal S}^{2,p}$.
\\ As a consequence of this and the fact that $(\int_0^{t}g_u du)_{t \in [0,T]}$ is in  ${\cal S}^{2,p}$ (since $g \in \H^{2}$),  we get that the martingale $N$ from the decomposition above belongs to ${\cal M}^2$. By the orthogonal decomposition property of martingales in ${\cal M}^2$ (cf. Lemma \ref{lemma1}), there exists a unique couple $(Z,M) \in \HB^2 \times  {\cal M}^{2,\perp} $ such that 
\begin{equation}
N_t =\int_0^t Z_s dW_s+M_t. \; \; \forall t \in [0,T]\;\;\;\mbox{a.s.}
\end{equation}
Thus, 
\begin{equation*}
Y_{\tau}=\xi_T+\int_{\tau}^T g(t)dt-\int_{\tau}^T  Z_t dW_t+M_{T^-}-M_{\t^-} +A_T-A_\tau+B_{T-} -B_{\tau-} \text{ a.s. for all }\tau\in\stopo.
\end{equation*}
	\item[(3)]\textbf{Minimality conditions}.
	Let $\tau \in \stopo$. It follows from (ii) Lemma \ref{Mertensvalue} that $\Delta B_{\t}= {\bf 1}_{\{ Y_\t= \xi_\t\}}\Delta B_{\t}$ a.s.\,
	which means that the process $B$ satisfies the minimality condition (v) of Definition \ref{defRBSDE}.\\ 
  On the other hand, we have from Lemma \ref{Mertensvalue} (iii) that
	$\Delta A_{\t}= {\bf 1}_{\{{ Y}_{\t-}= \, \xi_{\t^-}\}}\Delta
	A_{\t}$ a.s. We have also by Lemma \ref{lemc} that $\int_0^T 1_{\{Y_{t^-}>\xi_{t^-}\}}dA_t^c=0$. In other terms, $A$ satisfies the minimality condition (iv)  of Definition \ref{defRBSDE}.
	\item[(4)]$\boldsymbol{Y \geq \xi}$ up to an evanescent set. We use again the  aggregation equality of Lemma \ref{aggreg} combined with equation \eqref{ddeux-2} to obtain that,  $Y_\t=Y(\t)\geq \xi_\t$ for each $\t \in \stopo$. We obtain the desired result by an application of section theorem.
\end{description}
Collecting now all these properties yields that the quintuple $(Y,Z,M, A,B)$ is  a solution of the  following equation associated with driver  $g$ and the predictable obstacle $\xi$. That is
\begin{equation}\label{prince}
\left\{
\begin{array}{ll}
(i)&(Y,Z,M,A,B)\in {\cal S}^{2,p} \times \HB^2 \times {\cal M}^2 \times {\cal S}^{2,p}\times {\cal S}^{2,p}
\nonumber\\
(ii)& Y_\tau=\xi_T+\int_{\tau}^T g(t)dt-\int_{\tau}^T  Z_t dW_t+M_{T^-}-M_{\t^-} +A_T-A_\tau+B_{T-} -B_{\tau-} \text{ a.s. for all }\tau\in\stopo, \\ 
(iii)&  Y \geq \xi \text{ up to an evanescent set,}\\ \nonumber
(iv)&  A \text{ is a nondecreasing rightcontinuous predictable process with} \;A_0= 0, E(A_T)<\infty \text{ and such that } \\
&\int_0^T {\bf 1}_{\{Y_t > \xi_t\}} dA^c_t = 0 \text{ a.s. and }\;(Y_{\tau-}-\xi_{\tau-})(A^d_{\tau}-A^d_{\tau-})=0 \text{ a.s. for all (predictable) }\tau\in\stopo, \\   
(v)& B \text{ is a nondecreasing right-continuous adapted \textbf{ predictable} purely discontinuous} \\
&\text{ process with } B_{0-}= 0, E(B_T)<\infty  \text{ and such that }

(Y_{\tau}-\xi_{\tau})(B_{\tau}-B_{\tau-})=0 \text{ a.s. for all }\tau\in\stopo. 
\end{array}
\right.
\end{equation}
	\fproof
\subsection{Existence and uniqueness of the solution of Reflected BSDE associated with general driver}	
We are now going to prove that the RBSDE has a predictable solution in the general case of a general driver by using a fixed point argument with an appropriate mapping. We present also the link between RBSDEs in predictable framework and  predictable optimal stopping problem.  
	\begin{theorem}\label{existence general}
Let $\xi$ be a left limited predictable process in ${\cal{S}}^{2,p}$, let $g$ be a Lipschitz driver. The RBSDE associated with data $(g,\xi)$ has a unique solution $(Y,Z,M,A,B) \in {\cal S}^{2,p}  \times \HB^2\times {\cal M}^{2,p} \times {\cal S}^{2,p}\times {\cal S}^{2,p}$  and for each $S \in \stopo$, we have:
\begin{eqnarray}\label{char}
	Y_S= \esssup_{\t \in \stops} E\left[ \xi_\tau + \int_S^\t g(u,Y_u,  Z_u) du \mid \Fc_{S^-}\right].
\end{eqnarray}
	Moreover, the following  properties hold:
	\begin{description}
		\item[(ii)]	We have	$ Y_S=\xi_S\vee {}^pY^+_{S}$ a.s.
		for all $S \in \stopo$.
		\item[(iii)]	We have	$ Y_{S^-}=\xi_{S^-}\vee Y_{S}$ a.s. for all $S \in \stopo$.
	\item[(iv)] Moreover,
	For each $S \in  {\cal T}_{0,T}$ and for each $\alpha \in ]0,1[$, we set
	\[\tau^{\alpha}(S):= \inf\{t\geq S \,, \alpha Y_t+(1- \alpha)\int_0^tg(u,Y_u,Z_u)du\leq \xi_t \}.\]
	Then, 
	for each $\alpha  \in ]0,1[$ and for each $S \in \stopo$, 
	\begin{equation}\label{eq.mar}
	Y_S=E\left[-\int_0^{S} g(u,Y_u,  Z_u)du+\int_0^{\tau^\alpha(S)} Z_s dW_s+M_{\tau^\alpha(S)}-A_{\tau^\alpha(S)}-B_{\tau^\alpha(S)^{-}}|\mathcal{F}_{S^-}\right]\;\;\;\; \mbox{a.s.}
	\end{equation}
	Moreover,
	$A_{ \tau^\alpha(S)}=A_S$ and  $B_{\tau^\alpha(S)^-}=B_{S^-}$.
	\end{description}
\end{theorem}
\dproof
For each $\beta>0$, let $\mathbb{B}_\beta^{2,p}$ be the Banach space $ {\cal S}^{2,p}  \times \HB^2$, let $\psi$
a mapping defined from $\mathbb{B}_\beta^{2,p}$ onto itself as follows: for any $(U,V)\in \mathbb{B}_\beta^{2,p}$,  $(Y,Z)=\psi (U,V)$ is the unique element of  $\mathbb{B}_\beta^{2,p}$ such that $(Y,Z,M,A,B)$ solves the RBSDE associated with the driver $g(\omega,t):=g(t,\omega,U_t,V_t)$. The application is well defined by Theorem \ref{existence-linear}.

By using the a priori estimates from Lemma \ref{estimation} and similar computations as those from the proof of Theorem $3.4$ in \cite{MG}, we can prove that for a suitable choice of the parameter $\beta>0$, the mapping $\psi$ is a contraction from the Banach space into itself. By fixed point theorem , the mapping $\psi$ thus admits a unique fixed point, which corresponds to the unique solution to the RBSDE associated with $(g,\xi)$.\\
The assertions (ii), (iii) and (iv) follows from assertions (ii), (iii) and (iv) of Theorem \ref{existence-linear}, when the driver $g$ given by $g(\omega,t):=g(t,\omega,Y_t(\omega),Z_t(\omega))$.
\fproof	
\begin{Proposition}
Let $(g,\xi)$ a standard parameter associated with the unique solution $(Y,Z,M,A,B)$ of predictable RBSDE. Let $S \in \stopo$ and let $\beta$ be a nonnegative bounded $\mathcal{F}_{S^-}$-measurable random variable.
Then,  $(\beta Y,\beta Z,\beta M,\beta A,\beta B)$ is the unique solution on $[S,T]$ of the RBSDE with standard parameters $(\beta g,\beta \xi{\bf 1}_{]S,T]})$ 	
\end{Proposition}
\dproof
Let $Y^\beta$ be the first component  of the unique solution of RBSDE associated to  $(\beta g,\beta \xi{\bf 1}_{]S,T]})$. We have by equation \eqref{char}, for each $\tau \in \mathcal{T}_{S}^p $ 
 \begin{eqnarray*}
 	Y_\t^\beta= \esssup_{\theta \in \stops} E\left[ \beta \xi_\theta + \int_S^\theta \beta g(u,Y_u,  Z_u) du \mid \Fc_{\t^-}\right].
 \end{eqnarray*}
As a direct consequence of Proposition \ref{Propalpha}, $Y_\t^\beta=\beta Y_\t$  a.s. for all  $\tau \in \stopo$. If we rewrite the equation satisfied by $Y$ we can see that $(\beta Y,\alpha Z,\beta M,\beta A,\beta B)$ is the unique solution on $[S,T]$ of the RBSDE with standard parameters $(\beta g,\beta \xi{\bf 1}_{]S,T]})$. 

\begin{remark}
Let $S \in \stopo$ and $C \in\mathcal{F}_{S^-}$.  If we take $\beta={\bf 1}_{C}$, then,  $({\bf 1}_{C} Y,{\bf 1}_{C} Z,{\bf 1}_{C} M,{\bf 1}_{C}A,{\bf 1}_{C} B)$ the unique solution on $[S,T]$ of the RBSDE with standard parameters $({\bf 1}_{C} g,{\bf 1}_{C}\xi{\bf 1}_{]S,T]})$. 
\end{remark}
As we will see in the proof of Lemma \ref{estimation} in the Appendix, the estimates  still valid for non reflected BSDEs. Hence, we can use the same arguments  used in Theorem \ref{existence general}, to derive the following existence and uniqueness theorem:
\begin{theorem}
Let $\xi \in L^2(\cf_{T^-})$, let $g$ be a Lipschitz driver. The BSDE from Definition \ref{BSDE} with data $(g,\xi)$ has a unique solution.
\end{theorem}
\section{Optimal stopping with non linear predictable g expectation}\label{section-optimal}
 Let $(\xi_t)$ be a ladlag predictable process such that $\xi \in \mathcal{S}^{2,p}$, called \emph{reward}, modelling an agent's dynamic financial position. The agent's risk can be assessed by a dynamic risk measure induced by a predictable BSDE  with a given Lipschitz driver $g$; the dynamic risk measure equal up to a minus sign, to the predictable g-conditional expectation of $\xi$. If we consider an agent who can choose a predictable stopping time in $\mathcal{T}_0^p$, when she decide to stop at $\t \in \mathcal{T}_0^p$, she receives the amount $\xi_\t$ where $\xi_\t$ is a non negative $\mathcal{F}_{\t^-}$-random variable, and the risk is assessed by $-\mathcal{E}^{p,g}_{S,\t}(\xi_\t)$. The agent's aim  is to choose his strategy  in such a way that the risk be minimal as possible.\\
We formulate the optimal stopping problem  at time $0$ by:
\begin{equation}
V^p(0):=-{\rm ess}\sup_{\t\in \stopo}\mathcal{E}^{p,g}_{0,\t}(\xi_\t).
\end{equation}
For $S\in \stopo$, the \emph{ predictable value function at time $S$} is defined by the random variable 
\begin{equation}\label{eq.vs}
V^p(S):={\rm ess}\sup_{\t\in \mathcal{T}_S^p}\mathcal{E}^{p,g}_{S,\t}(\xi_\t).
\end{equation}
As mentioned in the introduction, the above optimal stopping problem over the set of stopping times, has been studied in \cite{EQ96} in the case of a continuous reward process $\xi$ and a Brownian filtration, in \cite{QuenSul2}  in the case of a right-continuous pay-off $\xi$, and in \cite{MG} in the case of a reward process which is only right-upper-semicontinuous, and in \cite{Imkeller2} without  any regularity assumptions on $\xi$. However, the study of optimal stopping problem \eqref{eq.vs} in the predictable context, seems also interesting  since it  gives us more information in modelling compared to classical cases.

\vline

We now introduce the following new definition,  which can be seen as an extension of the notion of predictable strong supermartingale processes.
\begin{definition}[Predictable strong ${\cal E}^{f}$-\emph{supermartingale family}]
	A process $(U_t)$ is
	said to be a predictable strong ${\cal E}^{p,g}$-\emph{supermartingale process} (resp. a  predictable strong ${\cal E}^{p,g}$ martingale process), if for any 
	$\sigma, \mu$ $ \in$ $\stopo$ such that $\mu \geq \sigma$ a.s.,
	\begin{eqnarray*}
		{\cal E}^{p,g}_{\sigma,\mu} (U_{\mu})  \leq  U_{\sigma} \quad \,\mbox{a.s.}&& {\rm (resp.} \quad 	{\cal E}^{p,g}_{\sigma,\mu} (U_{\mu}) =  U_{\sigma}\quad \,\mbox{a.s.}).
	\end{eqnarray*}
\end{definition} 
\begin{remark}
	Let $S \in \stopo$, $\t \in {\cal T}_0$ such that $S \leq \t$ a.s. The process $U$ is said to be a predictable 
	strong ${\cal E}^{p,g}$-\emph{  supermartingale process on} $[S,\t]$ (resp. a  predictable strong ${\cal E}^{p,g}$ martingale process on $[S,\t]$ ), if for all 
$ \sigma$, $\mu \in \stopo$  such that $S \leq \sigma \leq \mu \leq \t $ a.s., we have 
\begin{eqnarray*}
{\cal E}^{p,g}_{\sigma,\mu} (U_{\mu})  \leq  U_{\sigma} \quad \,\mbox{a.s.}&& {\rm (resp.} \quad 	{\cal E}^{p,g}_{\sigma,\mu} (U_{\mu}) =  U_{\sigma }\quad \,\mbox{a.s.} ).
\end{eqnarray*}
	
\end{remark}
\begin{remark}\label{rem.dec}
1	Let $U$ be  a predictable strong $\mathcal{E}^{p,g}$-supermartingale process. Let $S \in \stopo$,  then the application $\tau \rightarrow \mathcal{E}^{p,g}_{S,\tau}(U_{\tau}) $ is decreasing. In fact, let $\t, \t'\in \mathcal{T}_S^p$ such that $\t  \leq \t'$ a.s. 
	By the consistency property and monotonicity property of $\mathcal{E}^{p,g}$ we obtain,
	
	$$\mathcal{E}^{p,g}_{S,\t'}(U_{\t'})=\mathcal{E}^{p,g}_{S,\t}(\mathcal{E}^{p,g}_{\t,\t'}(U_{\t'}))\leq \mathcal{E}^{p,g}_{S,\t}(U_{\t}) \;\;{\rm a.s.}$$	
	Thus the result.	
\end{remark}
\begin{corollary}\label{cormartingale}
	
	Let  $ S \in \stopo$ and $A \in \mathcal{F}_{S^-}$. If the process  $(U_{t})_{t \in [0,T]}$ is a  predictable strong $\mathcal{E}^{p,g}$-martingale process, then the process $(1_AU_t)$ is  predictable strong $\mathcal{E}^{p,g1_A}-$  martingale process on $ [S,T]$.  
\end{corollary}
\dproof
	Let $\tau^1\leq  \tau^2 \in \stops$. Since  $S \leq \tau^1$, we have $A \in \mathcal{F}_{{\tau^1-}} $. Put $\t_A^2:= \t^2{\bf 1}_A+T{\bf 1}_{A^c}$. We
	have $\t_A^2 \in \mathcal{T}_{S}^p$ a.s. By using the zero-one law of   $g$-expectation which still holds in our  predictable setting (see  \cite[Proposition $15$]{Pe04}),  and that $\t_A^2=\t^2$ a.s. on $A$ and  the fact that $A \in \mathcal{F}_{{\tau_1}^-}$ together with the ${\cal E}^{p,g}$ martingale property of the process $(U_{t})$, give
	$$\mathcal{E}^{p,g1_A}_{\tau^1,\tau^2}( U_{\tau_2}{\bf
		1}_A)=\ce^{p,g^{\tau^2{\bf
				1}_A}}_{\tau^1,T}(U_{\tau^2}{\bf 1}_A)=\ce^{p,g^{\tau^2_A}}_{\tau_1,T}(U_{\tau^2_A}{\bf
		1}_A)=\ce^{p,g}_{\tau_1,\tau^2_A}(U_{\tau^2_A}){\bf
		1}_A=U_{\tau^1}{\bf
				1}_A.$$
	Which concludes the proof.
\fproof
\begin{Proposition}\label{pr.mar}
	Let $g$ be a Lipschitz driver. Let $S,\;\t \in \mathcal{T}_0^p$ with $S \leq \t$ a.s. Let $(U_{t})_{t \in [0,T]}$ be a  predictable strong $\mathcal{E}^{p,g}$-supermartingale process. The following assertions are equivalent:
	\begin{enumerate}
		\item  $ U_S=\mathcal{E}^{p,g}_{S,\t}(U_{\t})$ a.s.
		\item $\mathcal{E}^{p,g}_{0,S}(U_S)=\mathcal{E}^{p,g}_{0,\t}(U_{\t})$  a.s.
		\item   $\mathcal{E}^{p,g}_{0,\t'}(U_{\t'}) =\mathcal{E}^{p,g}_{0,\t}(U_{\t})$ a.s. for all $\t' \in \stops$ such that $\t' \leq  \t$ a.s.
		\item  $\mathcal{E}^{p,g}_{S,\t'}(U_{\t'})=\mathcal{E}^{p,g}_{S,\t}(U_{\t})$ a.s. for all $\t' \in \mathcal{T}_S$ such that $\t' \leq  \t$ a.s. 
		\item The process  $U$ is a predictable strong $\mathcal{E}^{p,g}$-martingale process on $[S,\t]$.
		
	\end{enumerate}
\end{Proposition}
\dproof
	$(1) \Rightarrow (2)$. The assumption  $U_S=\mathcal{E}^{p,g}_{S,\t}(U_{\t})$ a.s., together with consistency property of predictable $g$-expectations yield
	
	$$\mathcal{E}^{p,g}_{0,S}(U_{S})=\mathcal{E}^{p,g}_{0,S}(\mathcal{E}^{p,g}_{S,\t}(U_{\t}))=\mathcal{E}^{p,g}_{0,\t}(U_{\t}) \quad
	\mbox{\rm a.s.}.$$
	Which is the desired result.\\
	Let $\t ' \in \mathcal{T}_S^p$ such that $\t' \leq \t$ a.s.\\
	$(2)\Rightarrow (3)$. We will show that $\mathcal{E}^{p,g}_{0,\t'}(U_{\t'})=\mathcal{E}^{p,g}_{0,\t}(U_{\t})$. 
	By  Remark \ref{rem.dec}, the application 
	$\sigma \rightarrow \mathcal{E}^{p,g}_{0,\sigma}(U_{\sigma})$ is decreasing  on $[0,\t]$. Hence,
	$$\mathcal{E}^{p,g}_{0,\t}(U_{\t})\leq \mathcal{E}^{p,g}_{0,\t'}(U_{\t'}) \leq \mathcal{E}^{p,g}_{0,S}(U_{S}) \quad
	\mbox{\rm a.s.}.$$
	Finally, owing to assumption $2$, the last inequalities becomes equalities, and this completes the proof. \\
		$(3)\Rightarrow (4)$.
	 Under the hypothesis, $\mathcal{E}^{p,g}_{0,\t}(U_{\t})=\mathcal{E}^{p,g}_{0,\t'}(U_{\t'})$ a.s. By the consistency property of predictable $g$-expectations, this equality can be expressed as: $$\mathcal{E}^{p,g}_{0,S}(\mathcal{E}^{p,g}_{S,\t}(U_{\t}))=\mathcal{E}^{p,g}_{0,S}(\mathcal{E}^{p,g}_{S,\t'}(U_{\t'})) \quad
	\mbox{\rm a.s.}.$$
	On the other hand, we know that by Remark \ref{rem.dec}, $\mathcal{E}^{p,g}_{S,\t}(U_{\t}) \leq \mathcal{E}^{p,g}_{S,\t'}(U_{\t'})$ a.s.
	By the strict  monotonicity property of predictable $g$-expectations, we conclude that  $\mathcal{E}^{p,g}_{S,\t}(U_{\t}) = \mathcal{E}^{p,g}_{S,\t'}(U_{\t'})$ a.s.\\
	$(4)\Rightarrow (5)$.
	If $(4)$ holds, by the consistency property of predictable $g$-expectations, 
	$$\mathcal{E}^{p,g}_{S,\t'}
	(U_{\t'})=\mathcal{E}^{p,g}_{S,\t'}(\mathcal{E}^{p,g}_{\t',\t}(
	U_{\t}))\quad
	\mbox{\rm a.s.}.$$
	But, $\mathcal{E}^{p,g}_{\t',\t}(U_{\t}) \leq U_{\t'} $ a.s, thus, the strict monotonicity of predictable $g$-expectations permits us to deduce that $\mathcal{E}^{p,g}_{\t',\t}(U_{\t}) = U_{\t'}$ a.s., which is the desired result.\\
	$(5)\Rightarrow (1)$. Trivial.
\fproof
\begin{Proposition}\label{prop3}
Let $(U_t) \in \mathcal{S}^{2,p}$ be a right limited (RL) predictable strong $\mathcal{E}^{p,g}$-supermartingale process.
	For each $S, \, \t$ $\in$ $\mathcal{T}_0^p$ such that $\tau >S$, one has 
	\begin{itemize}
		\item  $\mathcal{E}^{p,g}_{S,\tau}
		(U_\tau) \leq {}^pU_{S}^+ $.
	\end{itemize}
\end{Proposition}
	\dproof
		Let $S, \, \t \in \mathcal{T}_0^p$ such that $\tau >S$, $(S+\frac{1}{n} \wedge \tau)$ is a predictable stopping time for all $n \in  \mathbb{N}$. By  the consistency property of $	\mathcal{E}^{p,g}$, the $\mathcal{E}^{p,g}$-supermartingale property, we get
		\begin{eqnarray*}
			\mathcal{E}^{p,g}_{S,\tau}
			(U_\tau)
			&\leq& \limn \mathcal{E}^{p,g}_{S,S+\frac{1}{n}\wedge \tau}\left(\mathcal{E}^{p,g}_{S+\frac{1}{n}\wedge \tau,\tau}\left( U_\t\right)\right) \leq \limn \mathcal{E}^{p,g}_{S,S+\frac{1}{n}\wedge \tau}\left[U_{S+\frac{1}{n}\wedge \tau}\right].
\end{eqnarray*}
Where the last inequality is due to the monotonicity of $\mathcal{E}^{p,g}$.
By using the definition of the operator $\mathcal{E}^{p,g}$  and the fact that $(U_t) \in \mathcal{S}^{2,p}$ combined with the (RL) property, one can check that:
$$\limn \mathcal{E}^{p,g}_{S,S+\frac{1}{n}\wedge \tau}\left(U_{S+\frac{1}{n}\wedge \tau}\right)={}^pU_{S}^+.$$
This concludes the proof.
\fproof
\begin{Proposition}\label{prop4}
	Let $(U_t) \in \mathcal{S}^{2,p}$ be a predictable strong $\mathcal{E}^{p,g}-$	supermartingale.  For each $S\in \stopo$, $\t \in \stops$, one has,
	\[	\mathcal{E}^{p,g}_{S,\tau} (U_{ \tau})\,1_{\{\tau >S\}}\leq {}^pU_{S}^+ 1_{\{\tau >S\}} \]
\end{Proposition}
\dproof
First, put $A=\{\tau >S\}$. Let us define the random variable $\overline \tau_A$  by $\overline\t_A:= \t{\bf 1}_A+T{\bf 1}_{A^c}$. Note that $A \in \mathcal{F}_{S^-}$, thus $\overline \tau_A$ belongs to $\mathcal{T}_{S^+}^p$. This with some  properties of predictable $g$-conditional expectation, we get
\[	\mathcal{E}^{p,g}_{ S,\tau} (U_{ \tau})\,1_{\{\tau >S\}} =\ce^{p,g^{\overline \tau{\bf
			1}_A}}_{S,T}[U_{\overline \tau}{\bf 1}_A]=\ce^{p,g^{\overline\tau_A}}_{S,T}[U_{\overline\tau_A}{\bf
	1}_A]=\ce^{p,g}_{S,\overline\tau_A}[U_{\overline\tau_A}]{\bf
	1}_A \,\,\,\mbox{\rm a.s.}\]
Since $\overline \tau_A>S$, we obtain by Proposition \ref{prop3}, $\ce^{p,g}_{S,\overline\tau_A}[U_{\overline\tau_A}]  \leq {}^pU_{S}^+$. Hence,
\[	\mathcal{E}^{p,g}_{ S,\tau} (U_{ \tau})\,{\bf 1}_{\{\tau >S\}}\leq {}^pU_{S}^+ \,{\bf 1}_{\{\tau >S\}}.\]
\fproof
\begin{Proposition}\label{propV_p}
	Let $(U_t) \in \mathcal{S}^{2,p}$ be a right limited predictable strong $\mathcal{E}^{p,g}$-supermartingale. For each $S \in \stopo$, one has 
	\begin{itemize}
		\item ${}^p U_{S}^+\leq U_{S}.$ a.s.
	\end{itemize}
\end{Proposition}
\dproof
	Let $S \in \mathcal{T}^p_0$, let $(S_{n})_{n} \in \mathcal{T}^p_S$ such that $S_n \downarrow S$. By the $\mathcal{E}^{p,g}$-supermartingale property,  $\mathcal{E}^{p,g}_{S,S^n}(U_{S^n})\leq U_{S}$ a.s. By using the definition of the operator $\mathcal{E}^{p,g}$  and the fact that $(U_t) \in \mathcal{S}^{2,p}$ combined with the (RL) property, we obtain
	$\limn \mathcal{E}^{p,g}_{S,S^n}(U_{S^n})={}^pU_{S}^+$. Thus,
	$${}^pU_{S}^+\leq U_{S}\;\;\;\;\mbox{a.s.}.$$
\fproof
\begin{Proposition}
	Let $(U_t) \in \mathcal{S}^{2,p}$ be a right limited predictable strong $\mathcal{E}^{p,g}$-supermartingale process. Then the process  $(^{p}U^+)$
	is also a predictable strong $\mathcal{E}^{p,g}$-supermartingale process.	
\end{Proposition}
\dproof Let $\tau,S \in \mathcal{T}^p_0 $ such that $\tau \geq S$, we have
	\begin{eqnarray}
\mathcal{E}^{p,g}_{S,\tau}
({}^pU_{\t}^+)&=	\mathcal{E}^{p,g}_{S,\tau}
	({}^pU_{\t}^+){\bf 1}_{\{\tau =S\}}+\mathcal{E}^{p,g}_{S,\tau}
	({}^pU_{\t}^+){\bf 1}_{\{\tau >S\}}.
	\end{eqnarray} 
	Let us remark that $$\mathcal{E}^{p,g}_{S,\tau}
	({}^pU_{\t}^+){\bf 1}_{\{\tau =S\}}={}^pU_{S}^+ {\bf 1}_{\{\tau =S\}}.$$ Moreover,  $$\mathcal{E}^{p,g}_{S,\tau}
	({}^pU_{\t}^+){\bf 1}_{\{\tau >S\}}\leq \mathcal{E}^{p,g}_{S,\tau}
	(U_{\t}){\bf 1}_{\{\tau >S\}}.$$
	This  last inequality is a consequence of Proposition \ref{propV_p} and the monotonicity property of $\mathcal{E}^{p,g}$. Thus,
		\begin{eqnarray}\label{eqle}
		\mathcal{E}^{p,g}_{S,\tau}
		({}^pU_{\t}^+)\leq 
		{}^pU_{S}^+{\bf 1}_{\{\tau =S\}}+\mathcal{E}^{p,g}_{S,\tau}
		(U_{\t}) {\bf 1}_{\{\tau >S\}}
		\end{eqnarray} 
	By Proposition \ref{prop3}, we have $$\mathcal{E}^{p,g}_{S,\tau}
	(U_\tau){\bf 1}_{\{\tau >S\}} \leq {}^pU_{S}^+{\bf 1}_{\{\tau >S\}} .$$ This combined with inequality \eqref{eqle} give $$\mathcal{E}^{p,g}_{S,\tau}
	({}^pU_{\t}^+)\leq {}^pU_{S}^+ .$$
	Thus, the desired result.
\fproof
\begin{Proposition}
Let $\xi$ be a ladlag predictable process. Let $(U_t) \in \mathcal{S}^{2,p}$ be a predictable strong $\mathcal{E}^{p,g}$-supermartingale,  such that $U \geq \xi$ up to an evanescent set. Then, for each $\t,\;\tilde\t \in \stopo$, such that  $\t\geq \tilde \t$, one has
\[	\mathcal{E}^{p,g}_{\tilde \t,\tau} (\xi_{ \tau})\,1_{\{\tau >\tilde \tau\}}\leq {}^pU_{\tilde \t}^+  1_{\{\tau >\tilde \tau\}} \]
\end{Proposition}
\dproof
We denote $A=1_{\{\tau >\tilde \tau \}}$. Let us define the random variable $\overline \tau_A$  by $\overline\t_A:= \t{\bf 1}_A+T{\bf 1}_{A^c}$. Note that $A \in \mathcal{F}_{S^-}$, thus $\overline \tau$ belongs to $\mathcal{T}_{\tilde \tau^+}^p$. 
Hence,
\[	\mathcal{E}^{p,g}_{\tilde \t,\tau} (\xi_{ \tau})\,1_{\{\tau >\tilde \tau\}} =\ce^{f^{\overline \tau{\bf
			1}_A}}_{\tilde \tau,T}[\xi_{\overline \tau}{\bf 1}_A]=\ce^{p,g^{\overline\tau_A}}_{\tilde \t,T}[\xi_{\overline\tau_A}{\bf
	1}_A]\leq \ce^{p,g^{\overline\tau_A}}_{\tilde \t,T}[U_{\overline\tau_A}{\bf
	1}_A]  \,\,\,\mbox{\rm a.s.}\]
Thus, By Proposition \ref{prop4}
\[	\mathcal{E}^{p,g}_{\tilde \t,\tau} (\xi_{ \tau})\,1_{\{\tau >\tilde \tau\}} \leq \ce^{p,g}_{\tilde \t,\overline\tau_A}[U_{\overline\tau_A}] 1_{\{\tau >\tilde \tau\}}\leq {}^pU_{\tilde \t}^+  1_{\{\tau >\tilde \tau\}}\]
Consequently, we get the desired result.
\fproof

\begin{lemma}\label{lemma_epsilon_optimality}

	Let $g$ be a Lipschitz driver and $\xi$ be ladlag process in $\mathcal{S}^{p,2}$.
	Let 
	$(Y,Z,M,A,B)$ be the solution to the reflected BSDE with parameters $(g,\xi)$ as in Definition \ref{defRBSDE}. Then,
	For each $S \in  {\cal T}_{0}^p$ and for each $\alpha \in ]0,1[$, the process $Y$ is a predictable strong  $\mathcal{E}^{p,g}$-martingale on $[S,\tau^{\alpha}_S]$.
\end{lemma} 
\dproof
By using Remark \ref{remark final} in the Appendix, it is  sufficient to show that $Y$ is the solution of the BSDE associated with driver $g$ and terminal time  $\t^\alpha(S)$ and terminal condition $Y_{\t^\alpha(S)}$. Thus, it is sufficient to prove that for a.e. $\omega \in \Omega$, the map $t\rightarrow A_t(\omega)+B_{t^-}(\omega)$ is constant on the closed interval $[S(\omega),\t^\alpha_S(\omega)]$.

In fact, The process $A$ is increasing  and we have by last assertion of Theorem \ref{existence general}, that $A_S=A_{\t^\alpha(S)}$ a.s. Thus, for a.e. $\omega \in \Omega$, $t\rightarrow A_t(\omega)$ is constant on  $[S(\omega),\t^\alpha_S(\omega)]$.
\\
On the other hand, $B$ is a non-decreasing right-continuous predictable purely discontinuous. Moreover, we have by the last assertion of Theorem \ref{existence general}, that $B_{S^-}=B_{\t^\alpha(S)^-}$. Hence, $B_{S^-}=B_{S}$ a.s. This shows that for a.e. $\omega \in \Omega$, $t\rightarrow B_{t^-}(\omega)$ is constant on  $[S(\omega),\t^\alpha_S(\omega)[$. By left-continuity of almost every trajectory of process $(B_{t^-})$  we obtain that for a.e. $\omega \in \Omega$, $t\rightarrow B_{t^-}(\omega)$ is constant on  $[S(\omega),\t^\alpha_S(\omega)]$. 
This concludes the proof.
\fproof

We introduce the following definition.
\begin{definition} \label{defr} A progressive process $(\xi_t)$ (resp. integrable) is said to be {\em right-upper semicontinuous (right USC) along stopping times (resp. along stopping times in expectation  (right USCE))} if for all $\tau \in {\cal T}_0$ and for each non increasing sequence of stopping times $ (\tau_n)$ such that $\tau^n \downarrow \tau$ a.s.\,,
	\begin{equation}\label{usc}
	\xi_{\tau} \geq \limsup_{n\to \infty} \xi_{\tau_n} \quad \mbox{a.s.} \quad (\text{ resp. } E[\xi_{\tau}] \geq \limsup_{n\to \infty} E[\xi_{\tau_n}] ).
	\end{equation}
\end{definition}
\begin{lemma}\label{lemma_epsilon_optimality}
	Let $g$ be a Lipschitz driver and $\xi$ be a positive left-limited right USCE process in $\mathcal{S}^{2,p}$.
	Let 
	$(Y,Z,M,A,B)$ be the solution to the reflected BSDE with parameters $(g,\xi)$. Let $\alpha>0$ and  $S \in \stopo$. Let
	$\theta^{\alpha}(S)$ be defined by
	\begin{equation}\label{eq_tau_epsilon_S}
	\theta^{\alpha}(S):= \essinf \{ \t  \in \stops \colon \alpha Y_\t \leq \xi_\t \}.
	\end{equation}
Then, we have
	$$\alpha Y_{	\theta^{\alpha}(S)}\leq  \;\xi_{\theta^{\alpha}(S)}\quad   \rm{a.s.}\,$$
\end{lemma} 
\dproof
Let $S \in \stopo$ and $A \in {\cal{F}}_{\theta^{\alpha}(S)^-}$.
In order to simplify notation, we denote $\theta^{\alpha}(S)$ by $\theta^{\alpha}$.
By definition of $\theta^{\alpha}$, there exists a non-increasing sequence $(\theta^n)$ in $\stops$  verifying  $\displaystyle{\theta^\alpha= \lim_{n\to \infty} \downarrow \theta^n}$  such that, we have  for each $n$,
\begin{equation}\label{eq.vn}
\alpha Y_{\theta^n}  \leq \xi_{\theta^n}\,\,\,{\rm a.s.}
\end{equation}

We have also by Remark \ref{rem upper}, \begin{equation}\label{eq.v+}
\alpha\; {}^p Y_{\theta^{n}}^+\leq \alpha Y_{\theta^n}\leq  \!\xi_{\theta^n}\,\,\,{\rm a.s.}
\end{equation} 
We have by (ii)  Theorem \ref{existence general}, 
$$E(\alpha\; Y_{\theta^{\alpha}}1_{A} )= E( \alpha\; {}^pY_{\theta^{\alpha}}^+1_{A \cap\{ Y_{\theta^\alpha}> \!\xi_{\theta^\alpha}\}})+E(\alpha \xi_{\theta^{\alpha}}1_{A \cap\{ Y_{\theta^\alpha}=  \!\xi_{\theta^\alpha}\}}).$$
Let us consider the first term of the r.h.s. of the last equality. 
Since $A \cap\{ Y_{\theta^\alpha}> \!\xi_{\theta^\alpha}\} \in {\cal{F}}_{\theta^{\alpha-}}$, we get, $$E(\alpha\; {}^pY_{\theta^{\alpha}}^+1_{A \cap\{ Y_{\theta^\alpha}> \!\xi_{\theta^\alpha}\}})=E(\alpha\; Y_{\theta^{\alpha+}}1_{A \cap\{ Y_{\theta^\alpha}> \!\xi_{\theta^\alpha}\}}).$$ 
On the other hand, the process $Y$ is in ${\cal S}^{2,p}$, thus 
\begin{equation*}
E(\alpha\; Y_{\theta^{\alpha+}}1_{A \cap\{ Y_{\theta^\alpha}> \!\xi_{\theta^\alpha}\}})=\limn E(\alpha\; Y_{\theta^n}1_{A \cap\{ Y_{\theta^\alpha}> \!\xi_{\theta^\alpha}\}}).
\end{equation*}
By inequality \eqref{eq.v+}, we obtain 
\begin{equation}
E(\alpha\; Y_{\theta^{\alpha+}}1_{A \cap\{ Y_{\theta^\alpha}> \!\xi_{\theta^\alpha}\}})\leq \limsupn E( \xi_{\theta^n}1_{A \cap\{ Y_{\theta^\alpha}> \!\xi_{\theta^\alpha}\}}).
\end{equation}
Thus, 
\begin{eqnarray*}
E(\alpha\; Y_{\theta^{\alpha}}1_{A})
&\leq& \limsupn E( \xi_{\theta^n}1_{A \cap\{ Y_{\theta^\alpha}> \!\xi_{\theta^\alpha}\}})+E( \xi_{\theta^{\alpha}}1_{A \cap\{Y_{\theta^\alpha}=\!\xi_{\theta^\alpha}\}})\\
&\leq& \limsupn E( \xi_{\bar \theta^n})
\end{eqnarray*}
where $\bar \theta^n:=\theta^n 1_{A \cap\{ Y_{\t^\alpha}> \!\xi_{\theta^\alpha}\}}+\theta^\alpha 1_{A \cap\{Y_{\theta^\alpha}=\!\xi_{\theta^\alpha}\}}+\theta^\alpha 1_{A^c} $.  Note that $\bar \theta^n$ is a non-increasing sequence of stopping times which verifies $\theta^\alpha=\displaystyle\lim_{n\to \infty} \downarrow \bar \theta^n$. Hence, by the right upper semicontinuity in expectations of the obstacle $\xi$, we obtain
\begin{equation}
E(\alpha\; Y_{\theta^{\alpha}}1_{A})\leq E(\xi_{\theta^\alpha}1_{A}).
\end{equation}
This holds for each $A \in {\cal{F}}_{\theta^{\alpha-}}$.
Thus,
  $$\alpha Y_{\theta^\alpha} \leq  \!\xi_{\theta^\alpha}.\quad a.s.$$
  Which is the desired result.
\fproof
\section{Existence of predictable optimal stopping time}\label{Existence of predictable}

\begin{definition}
 Let $S \in \stopo$. 
	Let $S$ $\in \stopo$ and let $\t_{*} \in \stopo$. We say that
		$\t_{*}$ is $S$-optimal for the value function $V_p$, if
		\begin{equation}\label{So}
	V_p(S) = \mathcal{E}^{p,g}_{S,\t_{*}}[\xi_{\t_{*}} \, ] \quad
		\,\mbox{a.s.}
		\end{equation}
\end{definition}

Classically, to prove the existence of an optimal stopping time, we prove the existence of an $\varepsilon$-optimal stopping time. This method is due to Maingueneau \cite{Maingueneau} and has found numerous applications (see for e.g. \cite{EK} in the case of a linear expectation), we refer also to \cite{QuenSul},\cite{QuenSul2}, \cite{MG}, in the case of non linear expectation. However, in our predictable setting, the use of this method lead to some additional complexities as for example the  exhibition of conditions on the process $\xi$ ensuring the existence of the solutions in this framework (see El Karoui \cite{EK}). To deal with this problem, we suggest another method to prove the existence of an optimal predictable stopping time. First, let us introduce the following set:
$$\mathcal{N}^p_S=\{\tau \in \mathcal{T}_S^p, \;\;\mbox{such that}\;\;(Y_{\sigma}, \sigma \in \mathcal{T}_{[S,\tau]})\;\;\mbox{is  a predictable  strong  $\mathcal{E}^{p,g}$- martingale}\}.$$ 
A natural candidate of optimal predictable stopping time for $Y$ is the random variable $ \tilde \tau(S) $  defined by
$$ \tilde \tau(S):=\esssup \mathcal{N}^p_S.$$
\begin{lemma}\label{stable}
	For each $S \in \stopo$, the set $\mathcal{N}^p_S$ is stable by pairwise maximization.
\end{lemma}
\dproof
	Let $S \in \stopo$ and $\tau_1,\;\tau_2 \in \mathcal{N}_S^p$. First, we have $\tau_1 \vee \tau_2  \in \mathcal{T}_S^p $. Let us show that  $\tau_1 \vee \tau_2  \in \mathcal{N}^p_S $. This is equivalent to show that $(Y_{\tau}, \tau \in \mathcal{T}_{[S,\tau_1 \vee \tau_2]})$  is a  predictable strong $\mathcal{E}^{p,g}$-martingale. By the consistency of predictable conditional $g-$expectation,  we get  a.s.
	\begin{eqnarray*}
		\ce^{p,g}_{S,\tau^1 \vee \tau^2}(Y_{\tau^1 \vee \tau^2})&=\ce^{p,g}_{S,\tau_1 \wedge \tau_2}(\ce^{p,g}_{\tau_1 \wedge \tau_2,\tau^1 \vee \tau^2}(Y_{\tau^1 \vee \tau^2})).\\
	\end{eqnarray*}	
	Let $A=1_{\{\tau^2>\tau^1\}}$, by  a standard property of predictable conditional $g-$ expectation,  the last equality can be rewritten as: 
	\begin{eqnarray*}
	\ce^{p,g}_{S,\tau^1 \vee \tau^2}(Y_{\tau^1 \vee \tau^2})=	\ce^{p,g}_{S,\tau_1 \wedge \tau_2}(\ce^{p,g^{\tau^1 \vee \tau^2 1_A}}_{\tau_1 \wedge \tau_2,T}({\bf 1}_A\;Y_{\tau^1 \vee \tau^2}))+\ce^{p,g}_{S,\tau_1 \wedge \tau_2}(\ce^{p,g^{\tau^1 \vee \tau^2 1_{A^c}}}_{\tau_1 \wedge \tau_2,T}({\bf 1}_{A^c}\;Y_{\tau^1 \vee \tau^2})).
	\end{eqnarray*}
	We have $\tau^1 \vee \tau^2 =\tau^2$ a.s. on $A$, $\tau^1 \vee \tau^2 =\tau^1$ a.s. on $A^c$. Therefore,
	\begin{eqnarray}\label{eqgen}
	\ce^{p,g}_{S,\tau^1 \vee \tau^2}[Y_{\tau^1 \vee \tau^2}]
	&=	\ce^{p,g}_{S,\tau_1 \wedge \tau_2}[\ce^{p,g^{ \tau^2 1_A}}_{\tau_1 \wedge \tau_2,T}[{\bf 1}_A\;Y_{\tau^2}]+[\ce^{p,g^{\tau^1  1_{A^c}}}_{\tau_1 \wedge \tau_2,T}[{\bf 1}_{A^c}\;Y_{\tau^1}]]\\
	&=\ce^{p,g}_{S,\tau_1 \wedge \tau_2}[{\bf 1}_A\ce^{p,g}_{\tau^1 \wedge \tau_2,\t^2}[Y_{\tau^2}]+{\bf 1}_{A^c}\ce^{p,g}_{\tau^1 \wedge \tau_2,\t^1}[Y_{\tau^1}]]
	\end{eqnarray}	
	Since $\t_1,\;\tau_2 \in \mathcal{N}_S$ and $A \in \mathcal{F}_{(\tau_1 \wedge \tau_2)^- }$, we have by Corollary \ref{cormartingale}, $({\bf 1}_AY_{\tau}, \tau \in \mathcal{T}_{[\tau_1\wedge \tau_2 ,\tau_2]})$ is  predictable strong $\mathcal{E}^{p,g1_A}$-  martingale process and $({\bf 1}_{A^c}Y_{\tau}, \tau \in \mathcal{T}_{[\tau_1\wedge \tau_2 ,\tau_1]})$ is  predictable strong $\mathcal{E}^{p,g1_{A^c}}$-  martingale process.
	Therefore, we get 
	\begin{eqnarray*}
		\ce^{p,g}_{S,\tau_1 \wedge \tau_2}[{\bf 1}_A\;Y_{\tau_1 \wedge \tau_2}+{\bf 1}_{A^c}\;Y_{\tau_1 \wedge \tau_1}]=\ce^{p,g}_{S,\tau_1 \wedge \tau_2}[Y_{\tau_1 \wedge \tau_2}]=Y_{S},
	\end{eqnarray*}
	where the last equality is due to the fact that  $(Y_{\tau}, \tau \in \mathcal{T}_{[S,\tau_1\wedge \tau_2]})$ is predictable strong $\ce^{p,g}-$  martingale process.
	Thus, $\ce^{p,g}_{S,\tau^1 \vee \tau^2}E[Y_{\tau^1 \vee \tau^2}]=Y_{S}$, by Corollary \ref{pr.mar} this give the desired result.
\fproof
\begin{theorem}\label{martingale}
	$ \tilde \tau(S)$ is a predictable stopping time. Moreover, 
assume that $(\xi_t)$ is  ${\cal{S}}^{2,p}$ which is l.u.s.c along stopping times. Then, the process $Y$ is a predictable strong $\ce^{p,g}$-maringale on $[S,\tilde\tau(S)]$	.
\end{theorem}
\dproof
	Let $S \in \stopo$, by Lemma \ref{stable} the set $\mathcal{N}^p_S$ is stable by pairwise maximization. Thus, there exists an increasing sequence $\tau_n$ of predictable stopping times in  $\mathcal{N}^p_S$ such that:
	$$\tilde \tau(S)=\lim \uparrow \tau_n.$$ Which proves that $\tilde \t$ is a predictable stopping time.
	We have for each $n \in \mathbb{N}$, $\tau_n \in \mathcal{N}^p_S$. Thus, for each $n$, $(Y_{\tau}, \tau \in \mathcal{T}_{[S ,\tau_n]})$ is a predictable strong $\ce^{p,g}$-martingale process. On the other hand, $Y$ is left continuous since $\xi$ is l.u.s.c along stopping times. Hence, by Proposition \ref{pr.mar} and  continuity of BSDEs with respect terminal time and terminal condition (see  \cite{QuenSul2}) which still holds in our predictable setting, we get
	$$Y_{S}=\limn \ce^{p,g}_{S,\tau_n}[Y_{\tau_n}]= \ce^{p,g}_{S,\tilde \tau(S)}[Y_{\tilde \tau(S)}].$$
	We can conculde by using again Proposition \ref{pr.mar}.
	\fproof
	\begin{theorem}\label{theroremmax}
		Let $g$ be a predictable Lipschitz driver. Let $(\xi_t)$ be a left-limited process in ${\cal{S}}^{2,p}$ which we assume to be l.u.s.c along stopping times and verifying, $ {}^p\xi^+ \leq \xi$. Let $(Y,Z,M,A,B)$ be the solution to the reflected BSDE with parameters $(g,\xi)$. Let $S \in \stopo$. Then, $Y_{\tilde \t_{S}}=\xi_{\tilde \t_{S}}$. a.s.
	\end{theorem}
	\dproof
	We note that $Y_{\tilde \t} \geq \xi_{\tilde \t}$ a.s., since $Y$ is the first component of the solution to the RBSDE with barrier $\xi$. We will show that $Y_{\tilde \t} \leq \xi_{\tilde \t}$ a.s. Suppose by the way of contradiction that  $P( Y_{\tilde\t_S} > \xi_{\tilde\tau})>0$. We have by the Skorohod condition satisfied by $B$, $\Delta B_{\tilde \tau} =B_{\tilde\tau}-B_{\tilde\tau^{-}}=0$ on the set $\{ Y_{\tilde \t}> \xi_{\tilde\tau}\}$. We have also that  $\Delta B_{\tilde\tau}={}^p Y_{\tilde\tau}^+-Y_{\tilde\tau}$. Hence, ${}^p Y_{\tilde\tau}^+=Y_{\tilde\tau}$ on the set  $\{ Y_{\tilde\tau}> \xi_{\tilde\tau}\}$. By definition of $\tilde \t$ as the essential spremum of  $\mathcal{N}^p_S$, we have 
		$Y_{\tilde\tau^{+}}= \xi_{\tilde\tau^{+}}$. Thus,
		$$Y_{\tilde\tau}={}^p Y_{\tilde\tau}^+={}^p \xi_{\tilde\tau^+} \leq \xi_{\tilde \t} \;\;\;\;\;\;\mbox{on the set} \;\;\;\{ Y_{\tilde\tau}> \xi_{\tilde\tau}\}.$$
		Which is a contradiction.
		\fproof
		\begin{theorem}\label{optimal}
			Let $g$ be a predictable Lipschitz driver. Let $(\xi_t)$ be a left-limited process in ${\cal{S}}^{2,p}$ which we assume to be l.u.s.c along stopping times and verifying, $ {}^p\xi_+ \leq \xi$. Let $(Y,Z,M,A,B)$ be the solution to the reflected BSDE with parameters $(g,\xi)$. Let $S \in \stopo$. Then, $\tilde \t_S$  is $S$- optimal for problem \ref{eq.vs}, that is 
			$$Y_S=\esssup_{\t\in \mathcal{T}_{S^p}}\mathcal{E}^{p,g}_{S,\t}(\xi_\t)=\mathcal{E}^{p,g}_{S,\tilde\t_S}(\xi_{\tilde\t_S}).$$
		\end{theorem}
		\dproof
		By Lemma \ref{compref}, the process $Y$  is a  predictable strong $\ce^{p,g}$-supermartingale. Hence, for each $\t \in\stops$, we have 
		$$Y_S \geq \mathcal{E}^{p,g}_{S,\t}(Y_\t) \geq \mathcal{E}^{p,g}_{S,\t}(\xi_\t). $$
		By taking the supremum over $\t \in \stops$, we obtain
		$$Y_S \geq {\rm ess}\sup_{\t\in \mathcal{T}_S^p}\mathcal{E}^{p,g}_{S,\t}(\xi_\t).$$
		It remains to show that $Y_S \leq {\rm ess}\sup_{\t\in \mathcal{T}_S^p}\mathcal{E}^{p,g}_{S,\t}(\xi_\t).$ We have by Theorem \ref{martingale}, the process $Y$ is a  predictable strong $\ce^{p,g}$-maringale on $[S,\tilde\tau(S)]$. Thus,  $Y_S=\mathcal{E}^{p,g}_{S,\tilde\tau(S)}(Y_{\tilde\tau(S)})$. On the other hand,  $Y_{\tilde \t_{S}}=\xi_{\tilde \t_{S}}$. a.s.by Theorem \ref{theroremmax}. Thus, we obtain,
		$$Y_S=\mathcal{E}^{p,g}_{S,\tilde\tau(S)}(\xi_{\tilde\tau(S)}),$$
		which yields 
		$$Y_S\leq {\rm ess}\sup_{\t\in \mathcal{T}_S^p}\mathcal{E}^{p,g}_{S,\t}(\xi_\t).$$

\fproof

We now in position to provide necessary and sufficient conditions, for optimal stopping time, in terms of appropriate martingales. This represents the non linear analogous in case of Bellman optimality criterium  (c.f El Karoui \cite{EK} in the setup of processes or \cite{Kob} Kobylanski and Quenez in the case of admissible families).
\begin{Proposition}\label{prop.criterion}\emph{(Optimality criterion) }
	Let $g$ be a Lipschitz driver and $\xi$ be a predicable barrier.
	Let 
	$(Y,Z,M,A,B)$ be the solution to the reflected BSDE with parameters $(g,\xi)$. 
	Let $S$ $\in \stopo$ and let $\t_{*} \in \stopo$.
	The following three assertions are equivalent
	\begin{enumerate}
		\item[(a)]
		$\t_{*}$ is $S$-optimal for $Y$, that is 
		\begin{equation}\label{So}
		Y_S = \mathcal{E}^{p,g}_{S,\t_{*}}[\xi_{\t_{*}} \, ] \quad
		\,\mbox{a.s.}
		\end{equation}
		\item [(b)] The following equalities hold:  
		$Y_{\t_{*}} = \xi_{ \t_{*}}\quad\mbox{a.s.,} \quad  {\rm and} \quad\mathcal{E}^{p,g}_{0,S}[Y_S]= \mathcal{E}^{p,g}_{0,\t_*}[Y_{\t_{*}}].$
		\item [(c)] The following equality holds: $\mathcal{E}^{p,g}_{0,S}[Y_S^p]= \mathcal{E}^{p,g}_{0,\t_{*}}[ \xi_{ \t_{*}}].$
	\end{enumerate}
\end{Proposition}
\dproof 
$(a) \Rightarrow (b)$: Since the process  $Y$ is a predictable strong $\mathcal{E}^{p,g}$-supermartingale process greater than $\xi$ (see Lemma \ref{compref}), and by monotonicity of predictable $g-$conditional  expectation,  we have clearly 
\[Y_S \geq \mathcal{E}^{p,g}_{S,\t*}[Y_{\t_{*}}] \geq  \mathcal{E}^{p,g}_{S,\t*}[\xi_{\t_{*}}]\;\mbox{a.s. } \]
By hypothesis, equality (\ref{So}) holds, this ensures that $\mathcal{E}^{p,g}_{S,\t*}[Y_{\t_{*}}] = \mathcal{E}^{p,g}_{S,\t*}[\xi_{\t_{*}}] $ a.s. On the other hand, the inequality 
$Y_{\t_{*}} \geq \xi_{\t_{*}}$ holds a.s. by the definition of $Y$. The strict monotonicity of predictable $g$- conditional expectation permits us to deduce that $Y_{\t_{*}} = \xi_{\t_{*}}$ a.s.
Moreover, $Y_S = \mathcal{E}^{p,g}_{S,\t*}[Y_{\t_{*}}] $ a.s., this combined with the consistency property of predictable $g$-conditional expectation   give $\mathcal{E}^{p,g}_{0,S}[Y_{S}]=\mathcal{E}^{p,g}_{0,S}[\mathcal{E}^{p,g}_{S,\t*}[Y_{\t_{*}}]]= \mathcal{E}^{p,g}_{0,\t*}[Y_{\t_{*}}]$, Hence, $(b)$ is satisfied. 
\\
$(b) \Rightarrow (c)$: it's clear.
\\
$(c) \Rightarrow (a)$: if $(c)$ holds, then by the consistency property of predictable $g$- conditional expectation, we can write $$\mathcal{E}^{p,g}_{0,S}[Y_S]=\mathcal{E}^{p,g}_{0,S}[\mathcal{E}^{p,g}_{S,\t*}[\xi_{\t_{*}}]]\,\,\,{\rm a.s.}.$$ 
Since $Y_S \geq\mathcal{E}^{p,g}_{S,\t_{*}}[\xi_{\t_{*}} ] $ a.s., the strict monotonicity of predictable $g$-conditional expectation give $Y_S = \mathcal{E}^{p,g}_{S,\t_{*}}[\xi_{ \t_{*}}] $ a.s.. Hence, $(a)$ is satisfied.
\fproof

\section{Some additional results on the strong predictable Snell envelope: the linear case}\label{section snell}
Let $\xi$ be a predictable reward process.
In this section, we study some properties of the predictable value function, defined at each predictable stopping time $S$ by
\begin{eqnarray}
V_p(S):= \esssup_{\t \in \stops} E\left[ \xi_\tau  \mid \Fc_{S^-}\right].
\end{eqnarray}
As in the seminal work of Kobylanski and Quenez \cite{Kob}, we avoid the aggregation step as well as the use of Merten's decomposition for strong predictable processes. Moreover, we only make assumption $\sup_{\t \in \stopo}E[|\xi_\t|]< \infty$ which is weaker than the assumption $E[\sup_{\t \in \stopo}|\xi_\t|]< \infty$  required in \cite{Karoui}.
\begin{definition}
	A family of random variables $\{\phi(\tau),\; \tau \in \mathcal{T}_0\}$ is said to be a predictable admissible family if it satisfies the following conditions:
	\begin{enumerate}
		\item for all $\tau \in \mathcal{T}_0^p$, $\phi(\tau)$ is an $\mathcal{F}_{\tau ^-}$-measurable $ \bar{\mathbb{R}}^+$-valued random variable,
		\item for all $\tau, \tau'\in \mathcal{T}_0^p $, $\phi(\tau)=\phi(\tau')$ a.s. on $\{\tau=\tau'\}$.
	\end{enumerate}
\end{definition}
In \cite{Karoui}, the reward is given by a predictable  process $(\phi_t)$. In this case, the family of random variables defined by $\{\phi(\tau)=\phi_{\tau},\;\tau \in \mathcal{T}_0^p\}$ is admissible.
\begin{Proposition}\label{increas}
	Given any two arbitrary predictable stopping times  $S$ and $\theta$ such that  $\theta \in  \mathcal{T}_{S}^p$, the family 	$\{E[\xi_{\t}|\mathcal{F}_{S^-}] \, , \, \tau\in\mathcal{T}_{\theta}^p\;\; \}$ is closed under pairwise maximization. Furthermore, 
	there exists a sequence of predictable stopping times $(\tau^n)_{n \in \mathbb{N}}$
	with $\tau^n $ in $ \mathcal{T}_{\theta}^p$  such that the sequence $(
	E[\xi_{\tau^n}|\mathcal{F}_{S^-}])_{n \in \mathbb{N}}$ converges non-decreasingly to $\displaystyle{\rm ess}\sup_{\tau\in \mathcal{T}_{\theta}^p} E[\xi_{\tau}|\mathcal{F}_{S^-}]$.	
\end{Proposition}
\dproof
	For any predictable stopping times $\tau^1$ and $\tau^2$  in $\mathcal{T}_{\theta}^{p}$,
	write\\ $A:=\{\, E[\xi_{\tau^2}|\mathcal{F}_{S^-}]\leq
	E[\xi_{\tau^1}|\mathcal{F}_{S^-}]\,\}$ and set
	$$\tau^3:=\tau^1 {\bf 1} _A+\tau^2 {\bf 1} _{A^c}.$$ 
	The fact that  $A \in \mathcal{F}_{S^-} \subset \mathcal{F}_{(\tau^1 \wedge \tau^2)^-}=\mathcal{F}_{(\tau^1)^-} \cap \mathcal{F}_{( \tau^2)^-}$, implies that $A \in\mathcal{F}_{(\tau^1)^- }$ and $A \in\mathcal{F}_{(\tau^2)^-}$. Thus, $\tau^3\in\mathcal{T}_{\theta }^{p}$, it follows that:\\
	$$ {\bf 1} _AE[\xi_{\tau^3})|\mathcal{F}_{S^-}]= 
	E[ {\bf 1} _A \xi_{\tau^3}|\mathcal{F}_{S^-}]= E[ {\bf 1} _A \xi_{\tau^1}|\mathcal{F}_{S^-}]
	= {\bf 1} _AE[\xi_{\tau^1})|\mathcal{F}_{S^-}]\;\;\;\;\mbox{ a.s}.$$
	Similarly, we show that $$ {\bf 1} _{A^c}E[\xi_{\tau^3}|\mathcal{F}_{S^-}]={\bf 1} _{A^c}E[\xi_{\tau^2}|\mathcal{F}_{S^-}]\;\;\;\; \mbox{ a.s}.$$ 
	Consequently,
	$$E[\xi_{\tau^3}|\mathcal{F}_{S^-}]=
	E[\xi_{\tau^1}|\mathcal{F}_{S^-}]{\bf 1} _A+
	E[  \xi_{\tau^2}|\mathcal{F}_{S^-}]{\bf 1} _{A^c}=
	E[\xi_{\tau^1}|\mathcal{F}_{S^-}]\vee
	E[\xi_{\tau^2}|\mathcal{F}_{S^-}] \mbox{ a.s},$$
	which shows the stability under pairwise maximization. Thus, by a classical result on essential supremum (see e.g. Neveu \cite{Neveu}), there exists a  sequence of predictable stopping times  $(\tau^n)_{n \in \mathbb{N}} \in\mathcal{T}_{\theta}^{p} $ such that
	
	$$\displaystyle{\rm ess}\sup_{\tau\in \mathcal{T}_{\theta}^{p}} E[\xi_{\tau}|\mathcal{F}_{S^-}]=\sup_n E[\xi_{\tau^n}|\mathcal{F}_{S^-}]\quad  \mbox{ a.s}., $$
	by recurrence, we can define a new sequence of stopping times $(\tilde\tau^n)_{n \in \mathbb{N}} \in\mathcal{T}_{\theta}^{p}$ by 
	$\tilde\tau^1=\tau^1$, and $\tilde\tau^n$ from $(\tilde\tau^{n-1},\tau^n)$ in the same way as in the definition of $\tau^3$ by $(\tau^{1},\tau^2)$. Hence, we can see that $E[\xi_{\tilde\tau^n}|\mathcal{F}_{S^-}]$ converges increasingly to $\displaystyle{\rm ess}\sup_{\tau\in \mathcal{T}_{\theta^+}^p} E[\xi_{\tau}|\mathcal{F}_{S^-}]$. 
	The proof is thus complete. 
\fproof
\begin{Proposition}\label{P1.2a}{\em(Optimizing sequences for $V_p$)}
	There exists a sequence of predictable stopping times $(\tau^n)_{n \in \mathbb{N}}$
	with $\tau^n $ in $ \mathcal{T}_{S}^p$ ,  such that the sequence $(E[\xi_{\tau^n}|\mathcal{F}_{S^-}])_{n \in \mathbb{N}}$ is increasing and such that
	\[V_p(S)= \lim_{n \to \infty} \uparrow E[\xi_{\tau^n}|\mathcal{F}_{S^-}] \quad
	\mbox{\rm a.s.}\]
\end{Proposition}
\dproof
	The result follows immediately by taking $\theta=S$ in Proposition \ref{increas}.
\fproof
\begin{lemma}\label{lemm3}
	Let $S \in \stopo $ and $ \theta \in \stops$. Let $\alpha $ be a nonnegative bounded $\mathcal{F}_{\theta^-}$-measurable random variable.  We have,
	\begin{eqnarray}\label{eqBellman}
	E[\alpha V_p(\theta)|\mathcal{F}_{S^-}]=\esssup_{  \tau \in \mathcal{T}_{\theta}^p}E[\alpha \xi_{\tau}|\mathcal{F}_{S^-}],
	\end{eqnarray}
\end{lemma}

\dproof
	 Let $\tau \in \mathcal{T}_{\theta}^p$,  by iterating expectation and using that $\alpha $ is a nonnegative bounded $\mathcal{F}_{\theta^-}$-measurable random variable,  combined  with  $E[ \xi_{\t}|\mathcal{F}_{\theta^-}] \leq  V_p(\theta)$, we obtain
	$$E[\alpha \xi_{\t}|\mathcal{F}_{S^-}]=E[E[\alpha \xi_{\t}
	|\mathcal{F}_{\theta^-}]|\mathcal{F}_{S^-}]=E[\alpha E[ \xi_{\t}|\mathcal{F}_{\theta^-}]|\mathcal{F}_{S^-}]\leq E[\alpha V_p(\theta)|\mathcal{F}_{S^-}].$$
	By taking the essential supremum over
	$\tau \in \mathcal{T}_{\theta}^p$ in the  inequality, we get 
	$$\esssup_{ \tau \in \mathcal{T}_{\theta}^p}E[\alpha \xi_{\t}|\mathcal{F}_{S^-}]\leq E[\alpha V_p(\theta)|\mathcal{F}_{S^-}].$$
	It remains to prove the reverse inequality $" \leq"$. By Proposition \ref{P1.2a},
	there exists a sequence of predictable stopping times $(\tau^n)_{n\in \mathbb{N}}$
	with 	 $\tau^n $ in $ \mathcal{T}_{\theta}^p$  and such that 
	$$V_p(\theta)=\limn \uparrow E[\xi_{\t^n}|\mathcal{F}_{\theta^-}].$$
	
	Since $\alpha$ is $\mathcal{F}_{\theta^-}-$measurable, we obtain that  $\displaystyle \alpha V_p(\theta)= \limn\uparrow E[ \alpha \xi_{\t^n}|\mathcal{F}_{\theta^-}]$ a.s. Therefore, applying the monotone convergence theorem and the fact that $S \leq \theta$ a.s. we derive that: \\
	$$E[\alpha V_p(\theta)|\mathcal{F}_{S^-}]=\limn \uparrow E[\alpha \xi_{\t^n}|\mathcal{F}_{S^-}].$$
	Hence,
	$$E[\alpha V_p(\theta)|\mathcal{F}_{S^-}] \leq \esssup_{  \tau \in \mathcal{T}_{\theta}^p}E[\alpha \xi_{\t}|\mathcal{F}_{S^-}].$$
	This with the previous inequality leads to the desired result.
\fproof

Let $\xi$ be  predictable reward process. Let $S \in \stopo$, let $\alpha$ be a nonnegative bounded $\mathcal{F}_{S^-}$-measurable random variable. Let $(V^\alpha(\tau),\; \tau \in \stopo)$ be the  value function associated with the reward $(\alpha\xi_{\theta},\;\theta \in \mathcal{T}_{S}^p)$, defined for each $\tau \in \mathcal{T}_{S}^p $ by 
$$V^\alpha(\tau):=\esssup_{\theta \in \mathcal{T}_{\tau}^p}E[\alpha \xi_\theta|\mathcal{F}_{\tau^-}].$$
Now, we will state some interesting properties:
\begin{Proposition}\label{Propalpha}
	Let $\xi$ be a predictable reward process, $S \in \stopo$ and let $\alpha$ be a nonnegative bounded $\mathcal{F}_{S^-}$-measurable random variable. The  value function $(V_p(\tau),\; \tau \in \stops)$  satisfies the following equality:
	\begin{itemize}
		\item $V^{\alpha}(\tau)=\alpha V_p(\tau)\;\;$ a.s. for all  $\tau \in \stops$.
	\end{itemize}
\end{Proposition}
\dproof
	Let	$\tau \in \stops$ and $\theta \in \mathcal{T}_{\tau}^p$. By the definition of the essential supremum (see Neveu \cite{Neveu} ),	$\alpha E[\xi_{\theta}|\mathcal{F}_{\tau^-}]=E[\alpha \xi_{\theta}|\mathcal{F}_{\tau^-}] \leq V^{\alpha}(\tau)$.
	Thus, by the characterization of the essential suprmem, we have 
	$\alpha V_p(\tau) \leq  V^{\alpha}(\tau) $. By the same arguments we can show that $V^{\alpha}(\tau) \leq  \alpha V_p(\tau) $. This concludes the proof.
\fproof

Let $S \in \stopo$ and $A \in\mathcal{F}_{S^-}$.  If we take $\alpha={\bf 1}_{A}$, we denote $V^\alpha$ by  $V^A$. Thus, $V^A$  is the value function associated with the  reward $(\xi_{\tau}1_A,\; \tau \in \mathcal{T}_{S}^p)$, defined for each $\tau \in \mathcal{T}_{S}^p $ 
$$V^A(\tau):=\esssup_{  \theta \in \mathcal{T}_{\tau}^p}E[ \xi_{\theta}1_A|\mathcal{F}_{\tau^-}].$$
\begin{lemma}\label{local}
	Let $\xi$ be a predictable reward process. Let $ \tau,\;\tilde \tau\in \stopo$ and denote $A:=\{\tau=\tilde \tau \}$. Then
	\begin{itemize}
		
\item $ V^A(\tau)=V^A(\tilde\tau)$\;\;\;\mbox{a.s}
\end{itemize}	
\end{lemma}
\dproof
	For each $\theta \in \mathcal{T}_{\tau}^p$, put $\theta_A=\theta {\bf 1}_A+T{\bf 1}_{A^c}$.
	Since $\tau$ and $\tilde \tau$ are predictable stopping times, we have $A \in \mathcal{F}_{\tau^-} \cap \mathcal{F}_{\tilde\tau^-}$. Thus, $\theta_A $ is predictable, thus we get a.s. on $A$:
	\begin{equation*}
	\begin{aligned}
	E[\xi_{\theta_A}{\bf 1}_A| \mathcal{ F}_{\tau^-}]&={\bf
		1}_AE[\xi_{\theta}| \mathcal{F}_{\tau^-}]={\bf
		1}_AE[\xi_{\theta}| \mathcal{F}_{\tilde\tau^-}]=E[\xi_{\theta_A}{\bf 1}_A| \mathcal{ F}_{{\tilde\tau}^-}],
	\end{aligned}
	\end{equation*}
	Since $\theta_A \in \mathcal{T}_{\tilde\tau}^p $, we obtain :
	$$E[\xi_{\theta_A}{\bf 1}_A| \mathcal{ F}_{\tau^-}]\leq  V^{A}(\tilde\tau).$$
	By arbitrariness of  $\theta  \in \mathcal{T}_{\tau^+}$,  this implies that 
	$$  V^{A}(\tau)\leq  V^{A}(\tilde\tau).$$
	By interchanging the roles of $\tau$ and $\tilde{\tau}$, we get $ V^{A}(\tau)= V^{A}(\tilde\tau)$.
	\fproof
	
Now, we will state te following localization property:
\begin{corollary}\label{corlocal}
	Let $(\xi)$ be reward process, $S \in \stopo$ and let $A \in \mathcal{F}_{S^-}$-measurable random variable. The  value function $(V_p(\tau),\; \tau \in \stops)$  satisfies the following equality:
	\begin{itemize}
		\item $V^{A}(\tau)={\bf 1}_{A} V_p(\tau)\;\;$ a.s. for all  $\tau \in \stops$.
	\end{itemize}
\end{corollary}
\dproof
	The result is a direct application of the Proposition \ref{Propalpha}.
\fproof

\begin{remark}\label{rem.p}
	Let $S \in \stopo$. Note that if $A \in\mathcal{F}_{S^-}$, we can always decompose the family  $(V_p(\tau),\; \tau \in \stops)$ as the following:
	$$ V_p(\tau)=V^A(\tau)+V^{A^c}(\tau) \;\;\;\;\;\;\; \mbox{for all}\;\;\;\;\;\tau \in \stops.$$ 
\end{remark}
The equalities above are useful, it allows us to prove the admissibility of the value function $V_p$.
\begin{Proposition}\label{P1.Adm}\emph{(Admissibility of $V_p$ )}\\
	The family $V_p=(V_p(S), S\in \stopo)$ is admissible.
\end{Proposition}
\dproof 
	For each $S\in\stopo, \; V_p(S)$ is an $\cf_{S^-}$-measurable random variable, due to the definition of the essential supremum (cf. e.g.  \cite{Neveu}).\\
	Let us prove Property $2$ of the  definition of admissibility. Take $\tau$ and $\tilde \tau$ in $\stopo$. We set $A:=\{\tau=\tilde \tau \}$ and
	we show that $V_p(\tau)=V_p(\tilde \tau)$, $P$-a.s. on $A$. \\
	Thanks to Lemma \ref{local}, $V^{A}(\tau)=V^{A}(\tilde\tau)$ a.s.
	Let us remark that $A \in \mathcal{F}_{\tau^- \wedge \tilde\tau^- } $.
	By  Corollary \ref{corlocal}, we have $$V_p(\tau){\bf 1}_{A}=V_p(\tilde\tau){\bf 1}_{A} \;\;\;\mbox{a.s.}$$	
	Thus the desired result.
\fproof
\begin{definition}[\emph{Predictable supermartingale system}]
	An admissible family $U:=(U(\tau), \; \tau \in\stopo)$ is
	said to be a  \emph{predictable supermartingale system} (resp. a  \emph{ predictable martingale system}) if, for any 
	$\tau, \tau^{'}$ $ \in$ $\mathcal{T}_0^p$ such that $\tau^{'} \geq \tau$ a.s.,
	\begin{eqnarray*}
		E[U(\tau')|{\cal{F}}_{\tau^-}] \leq  U(\tau) \quad \,\mbox{a.s.}&& {\rm (resp.}, \quad 	E[U(\tau')|{\cal{F}}_{\tau^-}]  =  U(\tau) \quad \,\mbox{a.s.}).
	\end{eqnarray*}
\end{definition} 

\vspace{0.3 cm}
A progressive process $X=(X_t)_{t \in [0,T]}$  is called a predictable strong supermartingale if it is a supermartingale, such that the family $(X_\tau,\; \tau \in \stopo)$ is a predictable supermartingale system.

\begin{corollary}\label{mar.loc}
	Let  $ S \in \stopo$ and $A \in \mathcal{F}_{S^-}$. If the family $(V_p(\tau),\; \tau \in \stops)$ is a predictable martingale system, then the family $(V^A(\tau),\; \tau \in \stops)$ is also  a predictable martingale system.  
\end{corollary}
\dproof
	Let $\tau_1< \tau_2 \in \stops$. Since  $S \leq \tau_1$, we have $A \in \mathcal{F}_{{\tau_1}^-} $. By applying Corollary \ref{corlocal}, and by using the martingale property of the system $(V_p(\tau),\; \tau \in \stops)$, we get  
	$$E[V^A(\tau_2)|\mathcal{F}_{\tau_{1}-}]=E[V_p(\tau_2)1_A|\mathcal{F}_{\tau_1-}]=E[V_p(\tau_2)|\mathcal{F}_{{
			\tau_1-}}]1_A=V_p(\tau_1)1_A=V^A(\tau_1).$$
	This concludes the proof.
\fproof

\begin{lemma}\label{lemm4}
	\begin{itemize}
		\item The admissible families $\{V_p(\tau), \tau \in \stopo\}$  is predictable supermartingale system. 
		\item The value family $V_p$ is characterized as the  predictable Snell envelope system associated with the reward process $\xi$, that is, the smallest supermartingale system which is greater (a.s.) than $\xi$.
	\end{itemize}
\end{lemma}
\dproof
 Let $S \leq \tau \in \stopo$. 
	Applying Lemma \ref{lemm3}, equation  \eqref{eqBellman} holds when $\alpha =1$. Since  $S \leq  \tau$, we get
	$$E[ V_p(\tau)|\mathcal{F}_{S^-}]=\esssup_{ \theta \in \mathcal{T}_{\tau}}E[ \xi_\theta|\mathcal{F}_{S^-}] \leq \esssup_{ \theta \in \mathcal{T}_S}E[ \xi_\theta|\mathcal{F}_{S^-}]=V_p(S),$$
	which gives the supermartingale property of $V_p$. \\
	Let us prove the second assertion. Let   $\{V_p'(\tau),\tau\in \stopo\}$  be another  supermartingale system
	such that $V_p'(\tau)\geq \xi_{\tau}$ a.s. for all $\tau
	\in\mathcal{T}_S^p$. Thus we have
	$$E[\xi_{\tau}|\mathcal{F}_{S^-}]\leq E[V_p'(\tau)|\mathcal{F}_{S^-}] \leq V_p'(S) \quad a.s.$$
	for all $\tau
	\in\mathcal{T}_S^p$. Hence
	by taking the essential supremum over
	$\tau \in\stops$, and by using the definition of $V_p$ we find that
	$$V_p(S)=\esssup_{ \tau \in \stops}E[\xi_{\tau}|\mathcal{F}_{S^-}]\leq V_p'(S) \quad a.s. $$ 
	for all $S \in\mathcal{T}_0^p$. This gives the desired result. 
\fproof

\section{Appendix}

\begin{lemma}\label{pro}	
	
	\begin{description}
		\item[(i)]  There exists a ladlag predictable process $( V_t)_{t\in[0,T]}$ which aggregates the family
		$(V(S))_{S\in\stopo}$ (i.e. $ V_S= V_p(S)$ a.s. for all $S\in\stopo$).\\
		Moreover, the process $(V_t)_{t\in[0,T]}$ is characterized as the predictable Snell envelope associated with the process 
		$( \xi_t)_{t \in [0,T]}$, that is the
		smallest predictable supermartingale
		greater than or equal to the process $( \xi_t )_{t \in [0,T]}$.
		\item[(ii)]	We have	$ V_S=\xi_S\vee {}^pV^+_{S}$ a.s.
		for all $S \in \stopo$.
		\item[(iii)]	We have	$ V_{S^-}=\xi_{S^-}\vee V_{S}$ a.s.
		for all $S \in \stopo$.
	\end{description}
\end{lemma}
\dproof	The proof is given in Theorem $2$ \cite{Karoui}.
\fproof
\begin{remark}\label{rem1}
	Let us remark that For all $S \in \stopo$, 	$V_{S}-{}^pV^+_{S}= {\bf 1}_{\{{ V}_S= \xi_S\}}(  V_{S}-{}^pV^+_{S})$ a.s.\,
	this follows from (ii) in the above Lemma. 
\end{remark}
\begin{remark}\label{rem2}
	We have for all  $S \in \stopo$, 	$V_{S^-}-V_{S}= {\bf 1}_{\{{ V}_{S^-}= \xi_{S^-}\}}(  V_{S}-V_{S})$ a.s.\,
	This is a direct consequence of iii).

\end{remark}
\begin{lemma}\label{principal}
Let $\xi$ be a predictable reward process such that $\xi \in {\cal S}^{2,p}$
	\begin{description}
		\item[(i)]  The predictable value process $( V_t)_{t\in[0,T]}$ is in ${\cal S}^{2,p}$ and admits the following predictable Mertens decomposition:
		\begin{equation}\label{eqmert}
		V_{\tau}=V_{0}+N_{\tau^-}-A_{\tau}-B_{\tau^-}\;\;\; \mbox{for all}\;\tau\in\stopo.
		\end{equation}
		where $N$ is a square integrable martingale, $A$ is a nondecreasing
		right-continuous predictable process such that $ B_0= 0$,
		$E(B_T^2)<\infty$, and $B$ is a nondecreasing right-continuous
	predictable purely discontinuous process such that  $B_{0-}= 0$,
		$E(B_T^2)<\infty$.
		\item[(ii)]
		For each $\tau \in\stopo$, we have
		$\Delta B_{\t}= {\bf 1}_{\{ V_\t= \xi_\t\}}\Delta B_{\t}$ a.s.\,
		\item[(iii)]  For each predictable  $\t \in\stopo$, we have
		$\Delta A_{\t}= {\bf 1}_{\{{ Y}_{\t-}= \, \xi_{\t^-}\}}\Delta
		A_{\t}$ a.s.
		\item[(iv)] For each $S \in  {\cal T}_{0,T}$ and for each $\alpha \in ]0,1[$, we set
		\[\tau^{\alpha}(S):= \inf\{t\geq S \,, \alpha { V}_t(\omega)\leq \xi_t \}.\]
		then 
		for each $\alpha  \in ]0,1[$, 
		$B_{ \tau^\alpha(S)}=B_S$ and  $B_{\tau^\alpha(S)^-}=B_{S^-}$.	
	\end{description}	
\end{lemma}
\dproof
By Lemma \ref{pro} (i), the process $V$ is a strong predictable supermartingale. By using martingales inequalities one can verify that 
\begin{equation}\label{eq_estimate_Y_f_xi}
E[\esssup_{S\in\stopo}|V_S|^2]\leq c \vertiii{\xi}^2_{{\cal S}^{2,p}}.
\end{equation}
Hence,  the process $( V_t)_{t\in[0,T]}$ is in
$\mathcal{S}^{2,p}$ (a fortiori, of class $(\mathcal{D}^p)$, i.e. $\{V_\t;\; \t \in \stopo\}$ is uniformly integrable). Applying Mertens
decomposition for predictable strong supermartingales of class ($\mathcal{D}^p$) (see \cite{Meyer_cours})
gives the decomposition \eqref{eqmert}, where $N$ is a cadlag uniformly integrable martingale, $A$ is a nondecreasing right-continuous predictable process
such that $ A_0= 0$, $E(A_T)<\infty$, and $B$ is a nondecreasing
right-continuous predictable purely discontinuous process such that
$B_{0-}= 0$, $E(B_T)<\infty$. By applying the same arguments as in
the proof of Lemma 3.3 (step 3) in \cite{Imkeller2}, one can verify that  $A
\in\mathcal{S}^{2,p}$ and $B \in\mathcal{S}^{2,p}$.\\
ii) Let $\tau \in\stopo$. It follows from the 
equation \eqref{eqmert} that:
$V_{\tau}^+=V_{0}+N_{\tau}-A_{\tau}-B_{\tau},$
hence,
${}^pV_{\tau}^+=V_{0}+N_{\tau^-}-A_{\tau}-B_{\tau}.$ This implies that $V_\tau-{}^pV_{\tau}^+=-(B_{\tau}-B_{\tau^-}).$ We conclude by Remark \ref{rem1} that $\Delta B_{\t}= {\bf 1}_{\{ V_\t= \xi_\t\}}\Delta B_{\t}\; a.s.$
\\
iii) We have by The Mertens decomposition \eqref{eqmert}, $V_{\tau^-}=V_{0}+N_{\tau^-}-A_{\tau^-}-B_{\tau^-}$ thus, $V_{S^-}-V_{S}=\Delta A_S,$. 
This combined with Remark \ref{rem2},  give $\Delta A_{\t}= {\bf 1}_{\{{ V}_{\t-}= \, \xi_{\t^-}\}}\Delta
A_{\t}\; a.s.$\\
iv)To sketch the proof, we refer the reader to Lemma $4$ in \cite{Karoui}.
\fproof
\begin{remark}
Since $A_{\tau^\alpha(S)}=A_S $ and $B_{\tau^\alpha(S)^-}=B_{S^-}$, we get
$$ V(S)= E[M_{\tau^\alpha(S)}-A_{\tau^\alpha(S)}-B_{\tau^\alpha(S)^{-}}|\mathcal{F}_{S^-}]. $$
\end{remark}	

\begin{theorem}\label{thmpr}
Let X be an $ \mathbb{\overline R}$ valued ${\cal F}\times \mathbb{R}_+$ measurable process. There exists an $\mathbb{R}$ valued process called the predictable projection of $X$ and denote $^pX$, that is determined uniquely up to an evansescent set by the following two conditions:
\begin{itemize}
	\item[(i)] it is predictable,
	\item[(ii)] $(^pX)=E[X_\t|\mathcal{F}_{\t^-} ]$ on $\{\t < \infty\}$. 
\end{itemize}
for all predictable times $\t$
\end{theorem}
\begin{corollary}\label{cor-predictable}
	If $X$ is a local martingale, then $(^pX)=X_-$.
\end{corollary}
	

\begin{theorem}[Gal'chouk-Lenglart]\label{Thm_Ito}
Let $n\in\N$. Let $X$ be an $n$-dimensional optional semimartingale, i.e. $X=(X^1,\ldots, X^n)$ is an $n$-dimensional optional process with decomposition $X^k=X_0^k+M^k +A^k+B^k$, for all $k\in\{1,\ldots,n\}$, where 
	$M^k$ is a left continuous local martingale, $A^k$ is a right-continuous process of finite variation such that $A_0=0$, and $B^k$ is a left-continuous process of finite variation which is purely discontinuous and such that $B_0=0$. Let $F$ be a twice continuously differentiable function on $\R^n$. Then, almost surely, for all $t\geq 0$,
	\begin{equation*}
	\begin{aligned}
	F(X_t)&=F(X_0)+\sum_{k=1}^{n}\int_{]0,t]} D^k F(X_{s-})d(A^k)_s \\
	&+\frac 1 2 \sum_{k,l=1}^{n}\int_{]0,t]} D^kD^l F(X_{s-})d[X^{k},X^{l}]_s\\
	&+ \sum_{0<s\leq t}\left[ F(X_s)-F(X_{s-})-\sum_{k=1}^{n} D^k F(X_{s-}) \Delta X_s^k-\frac 1 2 \sum_{k,l=1}^{n}\int_{]0,t]} D^kD^l F(X_{s-})\Delta X_s^k \Delta X_s^l\right]\\
	&+\sum_{k=1}^{n}\int_{[0,t[} D^k F(X_{s})d(B^k+M^k)_{s+}\\
	&+\sum_{0\leq s<t}\left[ F(X_{s+})-F(X_{s})-\sum_{k=1}^{n} D^k F(X_{s}) \Delta_+ X_s^k -\frac 1 2 \sum_{k,l=1}^{n}\int_{]0,t]} D^kD^l F(X_{s-})\Delta_+ X_s^k \Delta_+ X_s^l\right],
	\end{aligned}
	\end{equation*}
	where $D^k$ denotes the differentiation operator with respect to the $k$-th coordinate.
\end{theorem}

\begin{corollary}\label{Cor_Ito}
	Let $Y$ be a one-dimensional semimartingale with decomposition $Y=Y_0+M+A+B$, where 
	$M$, $A$, and $B$ are as in the above theorem. Let $\beta>0$. Then, almost surely, for all $t\geq 0$,
	\begin{equation*}
\begin{aligned}
\e^{\beta t}Y_t^2&=Y_0^2+\int_{]0,t]}\beta\e^{\beta s} Y_{s}^2 ds+2\int_{]0,t]} \e^{\beta s}Y_{s-}d(A)_s\\
&+\frac 1 2\int_{]0,t]} 2\e^{\beta s}d<M^{c},M^{c}>_s\\
&+ \sum_{0<s\leq t}\e^{\beta s}(Y_s-Y_{s-})^2\\	&+\int_{[0,t[} 2\e^{\beta s}Y_{s}d(B+M)_{s+}
+\sum_{0\leq s<t}\e^{\beta s}(Y_{s+}-Y_s)^2.
\end{aligned}
\end{equation*}
\end{corollary}
\dproof
	It suffices to apply the change of variables formula from Theorem \ref{Thm_Ito} with $n=2$, $F(x,y)=xy^2$, $X^1_t=\e^{\beta t}$, and $X^2_t=Y_t$. 
	Indeed, by applying Theorem \ref{Thm_Ito} and by noting  that the local martingale part and the purely discontinuous part of $X^1$ are both equal to $0$, we obtain
	\begin{equation*}
	\begin{aligned}
	\e^{\beta t}Y_t^2&=Y_0^2+\int_{]0,t]}\beta\e^{\beta s} Y_{s}^2 ds+2\int_{]0,t]} \e^{\beta s}Y_{s-}dA_s\\
&+\frac 1 2\int_{]0,t]} 2\e^{\beta s}d<M^{c},M^{c}>_s\\
&+ \sum_{0<s\leq t}\e^{\beta s}\big(Y_s^2-(Y_{s-})^2-2Y_{s-}(Y_s-Y_{s-})\big)\\	&+\int_{[0,t[} 2\e^{\beta s}Y_{s}d(B+M)_{s+}
+\sum_{0\leq s<t}\e^{\beta s}\big((Y_{s+})^2-(Y_{s})^2-2Y_s(Y_{s+}-Y_{s})\big).
\end{aligned}
\end{equation*}
	The desired expression follows as $Y_s^2-(Y_{s-})^2-2Y_{s-}(Y_s-Y_{s-})=(Y_s-Y_{s-})^2$ and 
	$(Y_{s+})^2-(Y_{s})^2-2Y_s(Y_{s+}-Y_{s})=(Y_{s+}-Y_s)^2.$

\dproof
	\textbf{Proof of Lemma \ref{estimation}:}
	Let $\beta>0$ and  $\varepsilon>0$ be such that $\beta\geq\frac 1 {\varepsilon^2}$.
	We  set $\tilde Y:=Y^1-Y^2$, $\tilde Z:=Z^1-Z^2$, $\tilde A:=A^1-A^2$, $\tilde B:=B^1-B^2$, and $\tilde g(\omega, t):=g^1(\omega, t)-g^2(\omega, t)$.
	We note that $\tilde Y_T=\xi_T-\xi_T=0;$ moreover, 
	$$\tilde Y_\tau=\int_{\tau}^T \tilde g(t)dt-\int_{\tau}^T  \tilde Z_t dW_t -\tilde M_{T^-}+\tilde M_{\t^-}+\tilde A_T-\tilde A_\tau+\tilde B_{T-} -\tilde B_{\tau-} \text{ a.s. for all }\tau\in\stopo.$$
	Thus we see that $\tilde Y$ is an optional (strong)  semimartingale (in the vocabulary of \cite{Galchouk}) with decomposition 
	$$\tilde Y_\t= \tilde Y_0+ M_-+A+B$$ where
	$M_t:=\int_0^t \tilde Z_s dW_s+\tilde M_t$, $A_t:= -\int_0^t  \tilde g(s) ds- \tilde A_t$ and 
	$B_t:=-\tilde B_{t-}$.
	
	Applying Corollary \ref{Cor_Ito} to $\e^{\beta t}\tilde Y_t^2$ gives: almost surely, for all $t\in[0,T]$,
	\begin{equation*}
	\begin{aligned}
	\e^{\beta T} (\tilde Y_T)^2&=\e^{\beta t}\tilde Y_t^2+\int_{]t,T]}\beta\e^{\beta s} (\tilde Y_{s})^2 ds-2\int_{]t,T]} \e^{\beta s}\tilde Y_{s-}\tilde g(s) ds -2\int_{]t,T]} \e^{\beta s}\tilde Y_{s-}d\tilde A_s\\ &+\int_{]t,T]} \e^{\beta s}d<M^{c},M^{c}>_s+\int_{]t,T]} \e^{\beta s} \tilde Z_s^2 ds
	+2\int_{]t,T]} \e^{\beta s}\tilde Y_{s-}\tilde Z_s d W_s
	\\
	&+2\int_{[t,T[} \e^{\beta s}\tilde Y_{s} d \tilde M_s-2\int_{[t,T[} \e^{\beta s}\tilde Y_{s} d \tilde B_s \\
	&+ \sum_{t<s\leq T}\e^{\beta s}(\tilde Y_s-\tilde Y_{s-})^2+\sum_{t\leq s<T}\e^{\beta s}(Y_{s+}-Y_s)^2.
	\end{aligned}
	\end{equation*}
	Since  $\tilde Y_T=0$ and $<\tilde M^c, W>=0$, we get: almost surely, for all $t\in[0,T]$,
	\begin{equation}\label{esproof}
	\begin{aligned}
	\e^{\beta t}\tilde Y_t^2
	&+\int_{]t,T]} \e^{\beta s} \tilde Z_s^2 ds+ \int_{]t,T]} \e^{\beta s}  d<\tilde M^c,\tilde M^c >_s=  -\int_{]t,T]}\beta\e^{\beta s} (\tilde Y_{s})^2 ds+
	2\int_{]t,T]} \e^{\beta s}\tilde Y_{s}\tilde g(s) ds\\
	&+2\int_{]t,T]} \e^{\beta s}\tilde Y_{s-}d\tilde A_s +2\int_{]t,T]} \e^{\beta s}\tilde Y_{s}d\tilde B_s-2\int_{]t,T]} \e^{\beta s}\tilde Y_{s-}\tilde Z_s d W_s-2\int_{[t,T[} \e^{\beta s}\tilde Y_{s} d \tilde M_s \\
	&-\sum_{t<s\leq T}\e^{\beta s}(\tilde Y_s-\tilde Y_{s-})^2-\sum_{0\leq s<t}\e^{\beta s}(Y_{s+}-Y_s)^2.
	\end{aligned}
	\end{equation}
	Let us first consider the sum of the first and the second term on the r.h.s. of the above inequality \eqref{esproof}. By the same arguments as in \cite{MG} and since $\beta\geq\frac 1 {\varepsilon^2}$, we get: a.e. for all $t\in[0,T]$,
	\begin{equation*}
	\begin{aligned}
	-\int_{]t,T]}\beta\e^{\beta s} (\tilde Y_{s})^2 ds+2\int_{]t,T]}  \e^{\beta s}\tilde Y_{s}\tilde g(s) ds\leq\varepsilon^2 \int_{]t,T]}  \e^{\beta s}\tilde g^2(s) ds.
	\end{aligned}
	\end{equation*}

	Next, we show that the third term  and fourth term on the right-hand side of inequality \eqref{esproof} are non-positive. More precisely,  a.s. for all $t\in[0,T]
	$,  
	$$\int_{]t,T]} \e^{\beta s}\tilde Y_{s-}d\tilde A_s \leq 0\;\;\;\;\mbox{and}\;\;\;\;\int_{[t,T[} \e^{\beta s}\tilde Y_{s}d(\tilde B)_{s}\leq 0.$$
	We give the detailed proof for the second inequality (the arguments for the first are similar). Indeed, a.s. for all $t\in[0,T]$,
	$\int_{[t,T[} \e^{\beta s}\tilde Y_{s}d(\tilde B)_{s}=
	\sum_{t\leq s<T} \e^{\beta s}\tilde Y_{s} \Delta \tilde B_s.$
	Now, a.s. for all $s\in[0,T]$,
	\begin{eqnarray*}
		\tilde Y_{s} \Delta \tilde B_s&=&(Y_s^1-Y_s^2)\Delta B_s^1-(Y_s^1-Y_s^2)\Delta B_s^2.
	\end{eqnarray*}
	First, we will show that $(Y_{s}^1-Y_{s}^2))\Delta B_s^1\leq 0$. We have
	
	$$(Y_{s}^1-Y_{s}^2)\Delta B_s^1=((Y_s^1-\xi_s)+ (\xi_s-Y_s^2))\Delta B_s^1=(Y_s^1-\xi_s)\Delta B_s^1 +(\xi_s-Y_s^2))\Delta B_s^1$$ By the Skorohod condition on $B^1$, $(Y_s^1-\xi_s)\Delta B_s^1=0$.
	By using the non-decreasingness of (almost all trajectories of) $B^1$  , and the fact that $Y^2\geq \xi$ we get: a.s. for all $s\in[0,T]$,
	$(Y_{s}^1-Y_{s}^2)\Delta B_s^1\leq 0.$ By similar arguments, we can show that $(Y_{s}^2-Y_{s}^1)\Delta B_s^2\leq 0.$ Thus, $\int_{[t,T[} \e^{\beta s}\tilde Y_{s}d(\tilde B)_{s}\leq 0$. By applying these observations to equation \eqref{esproof}, we get a.e. for all $t\in[0,T]$,
	\begin{equation}\label{eq-estimation}
	\begin{aligned}
	\e^{\beta t}\tilde Y_t^2&+ \int_{]t,T]} \e^{\beta s} \tilde Z_s^2 ds
	+ \int_{]t,T]} \e^{\beta s}  d<\tilde M^c,\tilde M^c >\\
	&\leq    \varepsilon^2 \int_{]t,T]}  \e^{\beta s}\tilde g^2(s) ds -2\int_{]t,T]} \e^{\beta s}\tilde Y_{s-}\tilde Z_s d W_s-2\int_{[t,T[} \e^{\beta s}\tilde Y_{s} d \tilde M_s \\
	&-\sum_{t\leq s< T}\e^{\beta s}(\tilde Y_{s^+}-\tilde Y_{s})^2.
	\end{aligned}
	\end{equation}	
	On the other hand, by definition of $\tilde Y$ we have, a.s. for all $s\in[0,T[$, $\Delta^+ \tilde Y_s=\Delta \tilde M_s-\Delta \tilde B_s$
	Hence, a.s. for all $t\in[0,T]$,
	\begin{eqnarray*}
	-\sum_{t\leq s< T}\e^{\beta s}(\tilde Y_{s^+}-\tilde Y_{s})^2
	=-\sum_{t\leq s< T}\e^{\beta s}(\Delta \tilde Ms-\Delta \tilde Bs)^2\leq -\sum_{t\leq s< T}\e^{\beta s}(\Delta \tilde M_s)^2+2\sum_{t\leq s< T}\e^{\beta s}\Delta \tilde M_s\Delta \tilde B_s
	\end{eqnarray*}
Hence, we obtain a.e. for all $t\in[0,T]$,
	\begin{equation*}
	\begin{aligned}
	\e^{\beta t}\tilde Y_t^2&+ \int_{]t,T]} \e^{\beta s} \tilde Z_s^2 ds+\int_{]t,T]} \e^{\beta s}  d<\tilde M^c,\tilde M^c >+\sum_{t\leq s<T}\e^{\beta s}(\Delta \tilde Ms)^2\\
	&\leq \varepsilon^2 \int_{]t,T]}  \e^{\beta s}\tilde g^2(s) ds-2\int_{]t,T]} \e^{\beta s}\tilde Y_{s-}\tilde Z_s d W_s-2\int_{[t,T[} \e^{\beta s}\tilde Y_{s} d \tilde M_s \\
	&+2\sum_{t\leq s<T}\e^{\beta s}\Delta \tilde M_s \Delta \tilde B_s.
	\end{aligned}
	\end{equation*}
	We have by the definition of the bracket term, $[\tilde M,\tilde M]=<\tilde M^c,\tilde M^c >+\sum_{s \leq .} (\Delta \tilde M)^2$
	We get: a.s. for all $t\in[0,T]$,
	Thus, we get
	\begin{equation}\label{eq-es-proof}
	\begin{aligned}
	\e^{\beta t}(\tilde Y_t)^2+\int_{]t,T]} \e^{\beta s} \tilde Z_s^2 ds+\int_{[t,T[} \e^{\beta s} d[\tilde M,\tilde M]_s &\leq \varepsilon^2 \int_{]t,T]}  \e^{\beta s}\tilde g^2(s) ds +2 (\bar{M}_{t^-}-\bar{M}_{T^-}).\\
	&+2\sum_{t<s\leq T}\e^{\beta s}\Delta \tilde M_s\Delta \tilde B_s.
	\end{aligned}
	\end{equation}
	where $\bar{M}$ is defined by:
	$$\bar{M}_t= 2\int_{]0,t]} \e^{\beta s}\tilde Y_{s-}\tilde Z_s d W_s+2\int_{]0,t]} \e^{\beta s}\tilde Y_{s} d \tilde M_{s}.$$
	We can verify that the local martingale $\bar{M}$ is a true martingale, by using some classical arguments based on the use of Burkholder-Davis-Gundy inequalities,
	Now, we will show that $E\left[\sum_{t\leq s<T}\e^{\beta s}\Delta \tilde M_s \Delta \tilde B_s=0\right]$. We note that $\tilde M$ is an uniformly integrable, thus $E\left[\Delta \tilde M_\t|{\cal{F}}_{\t^-}\right]=0$, for each predictable stopping time $\t \in \stopo.$ Moreover, $\Delta \tilde B_\t$ is a predictable process, since $\tilde B$ is predictable. Therefore, $\Delta \tilde  B_\t E\left[\Delta \tilde M_\t |{\cal{F}}_{\t^-}\right]=0$ or each predictable stopping time $\t \in \stopo.$ Hence,  $E\left[\sum_{t\leq s<T}\e^{\beta s}\Delta \tilde M_s \Delta \tilde B_s=0\right]=0$.

	By applying \eqref{eq-es-proof}with $t=0$ and by taking expectation on both sides, we get that 
	$$	(\tilde Y_0)^2+\|Z\|^2_{\beta}+E\left[\int_{[t,T[} \e^{\beta s} d[\tilde M,\tilde M]_s \right]\leq \varepsilon^2 \|\tilde g\|^2_{\beta}$$
	We get 
\begin{eqnarray}\label{estimate-Z}
\|Z\|^2_{\beta}\leq\varepsilon^2 \|\tilde g\|^2_{\beta} \;\;\;\;\;\;\mbox{and}\;\;\;\;E\left[\int_{[t,T[} \e^{\beta s} d[\tilde M,\tilde M]_s\right]\leq \varepsilon^2 \|\tilde g\|^2_{\beta}.
\end{eqnarray}
	
	\begin{equation}\label{eq-es}
	\begin{aligned}
	\e^{\beta t}(\tilde Y_t)^2 &\leq \varepsilon^2  \|\tilde g\|^2_{\beta} +2 (\bar{M}_t-\bar{M}_T).
	\end{aligned}
	\end{equation}	
	We get from \eqref{eq-estimation}
	\begin{equation*}\
	\begin{aligned}
	\e^{\beta t}\tilde Y_t^2\leq    \varepsilon^2 \int_{]t,T]}  \e^{\beta s}\tilde f^2(s) ds -2\int_{]t,T]} \e^{\beta s}\tilde Y_{s-}\tilde Z_s d W_s-2\int_{[t,T[} \e^{\beta s}\tilde Y_{s} d \tilde M_s 
	\end{aligned}
	\end{equation*}	
	By using Chasles's relation for stochastic integrals and by taking the essential supremum over $\tau\in\stopo$ and then the expectation on both sides of the above inequality , we obtain   
	\begin{equation}\label{eq9_lemma_estimate} 
	E[\esssup_{\tau \in \stopo}\e^{\beta \tau}\tilde Y_\tau^2]\leq \varepsilon^2 \|\tilde g\|^2_\beta  
	+2E[\esssup_{\tau \in \stopo}|\int_0^\tau \e^{\beta s}\tilde Y_{s-}\tilde Z_s d W_s|]+2E[\esssup_{\tau \in \stopo}|\int_{[0,\t[} \e^{\beta s}\tilde Y_{s} d \tilde M_s|].
	\end{equation}
	
Let us consider the third term of right hand side of the last inequality. By Burkholder-Davis-Gundy inequalities (applied with $p=1$), we get 
	$$E\left [\esssup_{\tau \in \stopo}|\int_{[0,\t[} \e^{\beta s}\tilde Y_s d \tilde M_s|\right] \leq c E\left[\sqrt{\int_{[0,T[} \e^{2\beta s}\tilde Y_{s}^2 d [\tilde M_s]|}\right] $$
	This combined with the inequality $ab \leq \frac{1}{2}a^2+\frac{1}{2}b^2$ yield
	\begin{eqnarray}
	2 E\left [\esssup_{\tau \in \stopo}|\int_{[0,\t[} \e^{\beta s}\tilde Y_s d \tilde M_s|\right]&\leq& E\left[\sqrt{\frac{1}{2}\esssup_{\tau\in\stopo} \e^{\beta_\t}(\tilde Y_{\tau})^2} \sqrt{8c^2\int_{[0,T[} \e^{\beta s}d[\tilde M]_s}\right]\\
	&\leq& \frac{1}{4}\|\tilde Y\|^2_{\beta}+4c^2E\left[\int_{[0,T[} \e^{\beta s} d[\tilde M]_s \right],
	\end{eqnarray}
	By using similar arguments, we obtain
	\begin{equation}\label{eq10_lemma_estimate} 
	2E[\esssup_{\tau \in \stopo}|\int_0^\tau \e^{\beta s}\tilde Y_{s-}\tilde Z_s d W_s|]
	\leq \frac{1}{4}\|\tilde Y\|^2_{\beta}+4c^2\|\tilde Z\|^2_{\beta}.
	\end{equation}
	From this, together with \eqref{eq9_lemma_estimate}, we get
	

	\begin{eqnarray*}\label{eq9_lemma_estimate2} 
		E[\esssup_{\tau \in \stopo}\e^{\beta \tau}\tilde Y_\tau^2]\leq \varepsilon^2 \|\tilde g\|^2_\beta  
		&+\frac{1}{2}\|\tilde Y\|^2_{\beta}+4c^2\|\tilde Z\|^2_{\beta}+4c^2E\left[\int_{[0,T[} \e^{\beta s} d[\tilde M]_s \right].
	\end{eqnarray*}
	    
	By this inequality, combined with the estimates \eqref{estimate-Z}, we get the following estimation
	$$\|\tilde Y\|^2_\beta \leq 2\varepsilon^2(1+8c^2)  \|\tilde f\|^2_\beta.$$  
	
\fproof
\begin{lemma}\label{compref} 
	Let $g$ be a  Lipschitz driver. Let $A$  be a nondecreasing 
	right-continuous predictable process in ${\cal S}^{2,p}$ with $A_0=0$
	and let $B$ be a nondecreasing 
	right-continuous predictable purely discontinuous process in ${\cal S}^{p,2}$ with $B_{0-}=0$.\\
	Let $(Y,Z,M) \in {\cal S}^{p,2} \times \mathbb{H}^2\times {\cal M}^{2, \bot}$  satisfy \\
	$$ -d  Y_t  \displaystyle =  g(t,Y_{t}, Z_{t})dt + dA_t + d B_{t-}-  Z_t dW_t- dM_{t^-}, \quad 0\leq t\leq T.$$
Then the process $(Y_t)$ is a predictable strong ${\cal E}^{p,g
}$-supermartingale. 
\end{lemma}

The proof use some specific arguments which are are suitable to the predictable setting as in the proof of precedent lemma and simliar arguments as those used in the proof  in \cite{MG}.
\begin{remark}\label{remark final}
We note that a process $Y\in{\cal S}^{p,2}$ is a strong  ${\cal E}^{p,g
}$-martingale on $[S,\t]$ (where $S,\t$ are two predictable stopping times  such that $S \leq  \t$ a.s.) if and only if, on $[S,\t]$, $Y$ is indistinguishable from the solution to the  BSDE from definition \ref{BSDE} associated with driver $g$, terminal time $\t$ and terminal condition $Y_\t$.
\end{remark}


\end{spacing}

\end{document}